\documentclass{article}
\usepackage{latexsym}
\usepackage[all]{xy}
\usepackage{amsmath}

\newtheorem{thm}{Th\'eor\`eme}[section]
\newtheorem{prop}[thm]{Proposition}
\newtheorem{lem}[thm]{Lemme}
\newtheorem{df}[thm]{D\'efinition}

\newtheorem{cor}[thm]{Corollaire}
\newtheorem{hyp}[thm]{Hypoth\`ese}

\begin{document}

\title{Th\'eor\`emes de Riemann-Roch pour les Champs de Deligne-Mumford}\date{}
\author{B. TOEN\thanks{Laboratoire Emile Picard, Universit\'e Paul Sabatier 118, route de Narbonne 31062 Toulouse Cedex France.
 e-mail:toen@picard.ups-tlse.fr}}
\maketitle

\begin{abstract}
On se propose de d\'evelopper un cadre cohomologique pour les
champs de Deligne-Mumford, adapt\'e \`a des formules de
type Hirzebruch-Riemann-Roch. On d\'efinit \`a cet effet la
"cohomologie \`a coefficients dans les repr\'esentations", ainsi
qu'un caract\`ere de Chern, et on d\'emontre un th\'eor\`eme
de Grothendieck-Riemann-Roch pour la transformation de Riemann-Roch
associ\'ee.\\

\bf Mots cl\'es: \rm Champs de Deligne-Mumford, th\'eor\`eme de Riemann-Roch, th\'eories cohomologiques. 
\end{abstract}

\newpage
\tableofcontents
\newpage

\begin{section}{Introduction}
\hspace{5mm}

Lorsque l'on \'etudie les espaces de modules, il est peu fr\'equent de 
rencontrer des espaces de modules fins. Ceci entraine en particulier que
certaines constructions naturelles ( les objets universels par exemple ) 
n'existent pas. Mais, si l'on se place dans le cadre plus g\'en\'eral des 
champs alg\'ebriques, ces constructions deviennent possibles. Il reste alors 
\`a en \'etudier les propri\'et\'es. Pour cela, il semble naturel de chercher \`a 
savoir si les th\'eor\`emes classiques de la g\'eom\'etrie alg\'ebrique restent 
valables lorsque l'on passe des vari\'et\'es aux champs. Nous nous proposons 
ici, d'analyser le cas du th\'eor\`eme de Grothendieck-Riemann-Roch.\\

Ennoncer un tel th\'eor\`eme n\'ecessite une th\'eorie cohomologique recevant un 
caract\`ere de Chern, et nous verrons que la principale difficult\'e r\'eside dans 
sa d\'efinition. \\
\hspace*{5mm}  Sur ce sujet il existe actuellement de nombreuses d\'efinitions de groupes 
de Chow d'un champ alg\'ebrique ( \cite{mu,eg,g2,vi2} ). Pour toutes ces 
th\'eories la projection naturelle d'un champ ( au sens de 
Deligne et Mumford \cite{dm} ) sur son espace de modules
$$p : F \longrightarrow M$$
v\'erifie $p_{*} : A(F)_{\bf Q} \simeq A(M)_{\bf Q}$. Nous verrons ( 
voir la remarque suivant \ref{p3a} ) que 
cette propri\'et\'e interdit une formule d'Hirzebruch-Riemann-Roch \`a 
valeurs dans $A(F)_{\bf Q}$. \\
\hspace*{5mm} Notre point de d\'epart est donc de trouver une nouvelle d\'efinition 
de la cohomologie d'un champ.
Pour cela, nous proposons d'\'etudier les 
spectres de $G$-th\'eorie des champs alg\'ebriques, afin d'un d\'eduire leurs 
comportements cohomologiques ( leurs "motifs" ). La d\'efinition s'imposera 
alors d'elle m\^eme. On retrouve ainsi une id\'ee de Carlos Simpson, 
consitant \`a consid\'erer des cycles dont les coefficients sont des repr\'esentations.\\

Le texte est d\'ecoup\'e en trois parties et un appendice.\\

Dans le premier chapitre on rappelle quelques r\'esultats et notations 
concernant la cohomologie des pr\'efaisceaux en spectres, et les champs 
alg\'ebriques. \\

Le second chapitre est consacr\'e \`a l'\'etude des spectres de $K$-th\'eorie des 
cat\'egories des faisceaux coh\'erents et localement libres sur un champ 
alg\'ebrique. On s'int\'eressera particulierement \`a deux probl\`emes.\\
\hspace*{5mm} Le premier concerne la descente de la $G$-th\'eorie rationnelle. Nous 
d\'emontrons \`a ce sujet trois r\'esultats
\begin{itemize}
\item Le th\'eor\`eme de descente \'etale, qui permet de travailler localement 
sur les espaces de modules.
\item Un r\'esultat de descente homologique analogue \`a celui d\'emontr\'e dans
\cite{g3}. Il permet par exemple de d\'efinir des images directes en 
$G$-th\'eorie \'etale, et de montrer par exemple que la projection naturelle
sur l'espace de modules
$$p : F \longrightarrow M$$
induit une \'equivalence faible 
$$p_{*} : H(F_{et},G\otimes\bf Q\mit) \simeq G(M)\otimes \bf Q$$
\item Enfin nous d\'efinissons une notion d'enveloppe de Chow 
rationnelle, qui 
permettra par la suite de nous ramener au cas de gerbes triviales.
\end{itemize}
\hspace*{5mm}Le second probl\`eme est la description des alg\`ebres $G_{*}(F)$, lorsque $F$ 
est lisse. Nous d\'emontrerons un isomorphisme fonctoriel de $\bf 
Q\mit(\mu_{\infty})$-alg\`ebres
$$G_{*}(F)\otimes \bf Q\mit(\mu_{\infty})
 \simeq H^{-*}((I_{F})_{et},G\otimes \bf Q\mit(\mu_{\infty}))$$
o\`u $I_{F}$ est le champ des ramifications de $F$. Notons que cet 
isomorphisme est une g\'en\'eralisation de la description de la $K$-th\'eorie 
\'equivariante donn\'ee par A. Vistoli dans \cite{vi}.\\

A la suite de cet isomorphisme, nous d\'efinirons dans le troisi\`eme 
chapitre, la cohomologie \`a coefficients dans les repr\'esentations par
$$H^{\bullet}_{rep}(F,*):=H^{\bullet}((I_{F})_{et},\Gamma(*))$$
o\`u $\Gamma$ est une th\'eorie cohomologique au sens de Gillet \cite{g}. 
On d\'efinit alors un caract\`ere de Chern et une classe de Todd, qui v\'erifient la 
formule de Grothendieck-Riemann-Roch
$$Ch(f_{*}(x)).Td(F')=f_{*}(Ch(x).Td(F))$$
pour un morphisme propre $f : F \longrightarrow F'$ de champs alg\'ebriques 
lisses. Comme dans le cas des sch\'emas, on diposera aussi d'une extension 
au cas des champs avec singularit\'es.\\ 

Dans l'appendice on donne une preuve d'un th\'eor\`eme de descente pour la 
cohomologie \`a coefficients dans les pr\'efaisceaux simpliciaux. C'est un 
r\'esultat que nous utilisons implicitement tout au long du texte. \\
\\
\\
\\

\underline{\bf REMERCIEMENTS:} \rm Je remercie tr\`es sinc\`erement
J. Tapia et C. Simpson, qui ont encadr\'e ce travail, et qui m'ont accord\'e de nombreuses discussions 
ainsi que de pr\'ecieux conseils.\\

Une importante partie du pr\'esent travail \`a \'et\'e effectu\'ee durant 
un s\'ejour \`a l'institut Max Planck, que je tiens \`a remercier pour la 
qualit\'e de son accueil.\\

\end{section}

\newpage

\begin{section}{Conventions et Notations}
\hspace{5mm}

On notera $\Delta$ la cat\'egorie simpliciale standart, et 
$SEns$ celle des ensembles simpliciaux. A l'aide de 
la construction $X \mapsto Sing(\mid X\mid)$, on supposera qu'ils sont 
tous fibrants.\\ 
\hspace*{5mm} De m\^eme $Sp$ d\'esigne la cat\'egorie des spectres ( \cite[1-1]{j} ), 
et on supposera que pour tout spectre $E$, toutes les composantes $E_{[n]}$ 
sont des ensembles simpliciaux fibrants.\\ 

\begin{subsection}{Cohomologie G\'en\'eralis\'ee}

\begin{subsubsection}{Pr\'efaisceaux Simpliciaux}
\hspace{5mm}
Soit $C$ un site. On d\'esigne par $SPr(C)$ la cat\'egorie des 
pr\'efaisceaux sur $C$ \`a valeurs dans $SEns$. D'apr\`es \cite[$2-32$]{j}, 
c'est une cat\'egorie de mod\`eles ferm\'ee simpliciale. 
Les termes fibration, cofibration et \'equivalence faible feront 
r\'ef\'erence \`a cette structure. Si $F$ et $G$ sont deux
objets de $SPr(C)$, on notera $Hom_{s}(F,G)$ l'ensemble simplicial des
morphismes de $F$ vers $G$ dans $SPr(C)$. Le pr\'efaisceau simplicial 
constant est not\'e $*$.\\
\hspace*{5mm} On dira qu'un morphisme $F \rightarrow G$ est une r\'esolution 
injective, si c'est une cofibration triviale, et si $G$ est un pr\'efaisceau 
simplicial fibrant. Tout pr\'efaisceau simplicial admet une 
r\'esolution injective, et elle est essentiellement unique.
En g\'en\'eral nous noterons $F \rightarrow HF$
une telle r\'esolution. On note alors, $H(X,F):=HF(X)$, pour chaque $X 
\in C$, et $H(C,F):=Hom_{s}(*,F)$. Si de plus $(F,x)$ est point\'e, on 
pose $H^{-m}(C,F):=\pi_{m}(H(C,F),x)$. 
Nous dirons que $F$ est un pr\'efaisceau simplicial flasque, si pour 
chaque objet $X \in C$, le morphisme $F(X) \rightarrow HF(X)$ est 
une \'equivalence faible.\\
\hspace*{5mm} Soit $U \rightarrow X$ un morphisme couvrant dans $C$, le nerf de $U$ 
sur $X$ est l'objet simplicial de $C$ d\'efini par
$$\begin{array}{cccc}
\cal N\mit(U/X) : & \Delta^{op} & \rightarrow & C \\
                  & [p]    & \mapsto     & \underbrace{U\times_{X}U\times_{X} 
                  \dots \times_{X}U}_{p+1 \; fois}
\end{array}$$
les faces et d\'eg\'en\'er\'escences \'etant donn\'ees par
les projections et les diagonales. Pour chaque morphisme couvrant $U \rightarrow X$, et $F$ un 
pr\'efaisceau simplicial, on dispose alors d'un $\Delta$-diagramme
$$F(\cal N\mit(U/X)) : [p] \mapsto F(U\times_{X}U\times_{X} 
                  \dots \times_{X}U)$$
La cohomologie de $\Check{C}$ech de $F$ pour le recouvrement 
$U \rightarrow X$ est l'ensemble simplicial
$$\Check{H}(U/X,F):=Holim_{\Delta}F(\cal N\mit(U/X))$$
o\`u $Holim$ est le foncteur de limite directe homotopique d\'ecrit 
dans \cite[Chap. $XI$]{bk}.\\
\hspace*{5mm} Nous dirons alors que $F$ v\'erifie la propri\'et\'e de Mayer-Vietoris, 
si pour chaque objet $X \in C$, et chaque recouvrement $U \rightarrow 
X$, le morphisme naturel $F(X) \rightarrow \Check{H}(U/X,F)$ est 
une \'equivalence faible. Le th\'eor\`eme de descente cohomologique ( 
voir l'appendice \ref{desc} ) nous dit qu'un pr\'efaisceau fibrant v\'erifie cette 
propri\'et\'e, et donc un pr\'efaisceau flasque aussi.\\                 

\end{subsubsection}

\begin{subsubsection}{Pr\'efaisceaux en Spectres}
\hspace{5mm}
On notera $SpSpr(C)$ la cat\'egorie des pr\'efaisceaux en spectres 
sur un site $C$. D'apr\`es \cite[$2-34$]{j}, c'est une cat\'egorie de mod\`eles ferm\'ee.
On dispose donc de notions de fibration, cofibration et 
\'equivalence faible. Si $F$ et $G$ sont deux pr\'efaisceaux en 
spectres on note $Hom_{sp}(F,G)$ le spectre des morphismes de $F$ 
vers $G$. Comme pr\'ec\'edemment, les  r\'esolutions injectives seront not\'ees
$F \rightarrow HF$, et $H(X,F):=HF(X)$. Si $F$ est un pr\'efaisceau 
simplicial, on note $SF$ son spectre de suspension. On d\'efinit 
$\Sigma:=S*$. Si $F$ est un pr\'efaisceau en spectres, on pose 
$H(C,F):=Hom_{sp}(\Sigma,F)$, et $H^{-m}(C,F):=\pi_{m}(H(C,F))$.
On dispose aussi de notions de cohomologie de $\Check{C}$ech, de 
pr\'efaisceaux en spectres flasques, et du th\'eor\`eme de descente.\\
\hspace*{5mm} Soit $f : F \rightarrow G$ un morphisme de pr\'efaisceaux en 
spectres. La fibre homotopique de $f$ est d\'efinie par
$$Fib(f)_{[n]}:=Holim \left(\xymatrix{
F_{[n]} \ar[r]^{f} & G_{[n]} \\
  & \ast \ar[u]} \right)$$
La cofibre homotopique est 
$$Cof(f)_{[n]}:=Hocolim \left( \xymatrix{
F_{[n]} \ar[d] \ar[r]^{f} & G_{[n]} \\
  \ast & } \right)$$ 
Une suite exacte est la donn\'ee de morphismes de pr\'efaisceaux en spectres
$$\xymatrix{
E \ar[r]^{f} & F \ar[r]^{g} & G }$$
et d'un diagramme commutatif
$$\xymatrix{
E \ar[r]^{f} \ar[d]_{u} & F \ar[r]^{g} \ar[d]_{id} & G \ar[d]_{id} \\
Fib(g) \ar[r]^{j} & F \ar[r]^{g} & G \\}$$
o\`u $j$ est le morphisme 
naturel, et $u$ une \'equivalence faible.
D'apr\`es \cite[$4-1$, $4-4$]{j}, c'est aussi \'equivalent \`a la donn\'ee d'un diagramme commutatif
$$\xymatrix{
E \ar[r]^{f} & F \ar[r]^{g}  & G  \\
E \ar[u]^{id} \ar[r]^{f} & F \ar[u]^{id} \ar[r]^{p} & Cof(f) \ar[u]^{v}\\}$$
avec $p$ le morphisme naturel et $v$ une \'equivalence faible.\\
\hspace*{5mm} Un morphisme de suites exactes est un diagramme commutatif
$$\xymatrix{
E \ar[r]^{f} \ar[d] & F \ar[r]^{g} \ar[d] & G \ar[d] \\
E' \ar[r]^{f'} & F' \ar[r]^{g'} & G' \\}$$
compatible avec les diagrammes ci-dessus.\\
\hspace*{5mm} Toute suite exacte 
$$\xymatrix{
E \ar[r]^{f} & F \ar[r]^{g} & G }$$
donne une suite exacte longue de cohomologie
$$\xymatrix{
 H^{-m}(C,F) \ar[r]^{g_{*}} & H^{-m}(C,G) \ar[r]^{\delta} & H^{-m+1}(C,E) 
\ar[r]^{f_{*}} & H^{-m+1}(C,F) }$$
Un morphisme de suites exactes induit un morphisme de suites exactes 
longues. Remarquons enfin que le foncteur $Holim$ ( resp. 
$Hocolim$ ) commute avec la formation des fibres homotopiques ( resp. cofibres 
homotopiques ), et transforment donc de fa\c{c}on naturelle les suites 
exactes en suites exactes.\\ 

\end{subsubsection}

\end{subsection}

\begin{subsection}{Champs Alg\'ebriques}
\hspace{5mm}

Les r\'ef\'erences pour ce paragraphe sont \cite{dm,lm}.\\
 
Soit $p : \cal C\mit \rightarrow C$ une cat\'egorie fibr\'ee en groupo\`edes sur un site 
$C$. On lui associe un pr\'efaisceau en groupo\`edes sur $C$ 
d\'efini par
$$\begin{array}{cccc}
F_{\cal C\mit} : & C & \rightarrow & Groupoides \\
                 & U & \mapsto & Hom_{C}(U, \cal C\mit) 
\end{array}$$ 
avec $U$ consid\'er\'e comme cat\'egorie fibr\'ee repr\'esent\'ee par $U$, et
$Hom_{C}$ le groupo\`ede des morphismes de cat\'egories 
fibr\'ees sur $C$. En associant \`a chaque groupo\`ede son classifiant, on 
obtient un pr\'efaisceau simplicial $F_{\cal C\mit}$ sur $C$. Avec ces 
notations, $\cal C\mit$ est un champ si $F_{\cal C\mit}$ est flasque 
comme pr\'efaisceau simplicial. \\
\hspace*{5mm} De cette fa\c{c}on on a une \'equivalence de la $2$-cat\'egorie des champs 
sur $C$ avec celle des pr\'efaisceaux simpliciaux flasques $1$-tronqu\'es 
et morphismes flexibles ( \cite{s} ).\\ 
\hspace*{5mm} Si $F \rightarrow H$ et $G \rightarrow H$ sont deux morphismes de 
champs, on note $F\times_{H}G$ le produit fibr\'e homotopique. Si les 
morphismes vers $H$ sont clairs, on notera aussi $F\times G$. \\
\hspace*{5mm} Soit $f : F \rightarrow G$ un morphisme de champ, et $U \rightarrow 
G$ un morphisme avec $U \in C$. La fibre de $f$ au-dessus de $U$ est 
not\'ee $f^{-1}(U):=F\times_{G}U$. Rappelons que $f$ est repr\'esentable
si toutes les fibres $f^{-1}(U)$ sont repr\'esentables
par des objets de $C$.\\

Supposons maintenant que $C=(Sch/k)_{et}$, le site des sch\'emas s\'epar\'es 
de type fini sur un corps $k$. Un sch\'ema sera un objet de $C$.
Un champ $F$ sera alg\'ebrique si
d'une part le morphisme diagonal 
$$\Delta : F \rightarrow F\times F$$
est repr\'esentable et fini, et de plus s'il existe un sch\'ema $X$ et un 
morphisme \'etale et surjectif $X \rightarrow F$.
Le champ des ramifications d'un champ alg\'ebrique $F$ est 
$$I_{F}:=F\times_{F\times F}F$$
Il est muni d'un morphisme $\pi : I_{F} \rightarrow F$ qui
est ( repr\'esentable ) fini et non-ramifi\'e. On dira que $F$ est une gerbe si
$\pi$ est \'etale.\\

Soit $X$ une vari\'et\'e et $H$ un groupe fini op\'erant sur $X$. Le champ
quotient de $X$ par $H$ est not\'e $[X/H]$. Ses sections au-dessus d'un
sch\'ema $Y$ sont les $H$-torseurs $P$ sur $Y$ munis de morphismes 
$H$-\'equivariants vers $X$. Si $F=[X/H]$, alors on a une \'equivalence 
canonique $I_{F}\simeq[\Hat{X}/H]$, avec 
$\Hat{X}=\{(x,h) \in X\times H / h.x=x \}$. Et donc
$$I_{F}\simeq \coprod_{h \in c(H)}[X^{h}/Z_{h}]$$
avec $c(H)$ les classes de conjugaisons de $H$, et $Z_{h}$ le 
centralisateur de $h$ dans $H$.\\

Un espace de modules pour $F$, est un espace alg\'ebrique $M$ muni 
d'un morphisme propre $p : F \rightarrow M$ qui est universel vers
les espaces alg\'ebriques, et tel que $p : \pi_{0}F(Spec K) \simeq 
M(Spec K)$ soit une bijection pour tout corps alg\'ebriquement clos 
$K$. D'apr\`es \cite{km}, tout champ alg\'ebrique poss\`ede un espace de 
modules. De plus, il existe un recouvrement \'etale $U \rightarrow M$, un sch\'ema 
$V$, un faisceau constant de groupes finis $H$ op\'erant sur $V$, et 
une \'equivalence $F_{U}:=p^{-1}(U)\simeq [V/H]$. Le champ $F$ est 
une gerbe si et seulement si on peut prendre $H$ qui op\`ere trivialement 
( on dit alors que " $H$ est le groupe de la gerbe $F$ " ). En 
particulier, si $F$ est r\'eduit, il existe un sous-champ ouvert dense qui 
est une gerbe. Les sous-champs ouverts ( resp. ferm\'es r\'eduits ) de $F$ sont 
tous de la forme 
$F':=F\times_{M}M' \hookrightarrow F$ ( resp. $F':=(F\times_{M}M')_{red} 
\hookrightarrow F$ ), 
 avec $M' \hookrightarrow M$ une immersion 
ouverte ( resp. ferm\'ee ). On notera $Codim(F',F):=Codim(M',M)$.
On appellera point de $F$ un point de $M$. La gerbe r\'esiduelle d'un 
point $Spec k(x) \rightarrow M$, est la gerbe sur $Spec k(x)$, 
$\widetilde{x}=(F\times_{M}Spec k(x))_{red}$. L'ordre d'inertie de $F$ en $x$ 
est l'ordre du groupe 
de la gerbe r\'esiduelle $\widetilde{x}$. \\
\hspace*{5mm} Soit $M' \rightarrow M$ est un morphisme d'espaces alg\'ebriques, $M$ 
l'espace de modules de $F$, et $F'=F\times_{M}
M'$. Alors le morphisme naturel de l'espace de modules $M''$ de $F'$
$$M'' \rightarrow M'$$
est radiciel ( \cite[$2.1$]{vi} ). Comme nous nous interresserons aux types 
d'homologie rationnelle, on pourra indentifier $M''$ et $M'$ par ce 
morphisme.\\
  
Le site \'etale d'un champ alg\'ebrique $F$ sera not\'e $F_{et}$ ( 
\cite[$4.10$]{dm} ). Il est muni d'un faisceau d'anneaux $\cal O\mit_{F}$. 
La cat\'egorie des faisceaux de $\cal O\mit_{F}$-modules coh\'erents sera
not\'ee $Coh(F)$, et celle des faisceaux de $\cal O\mit_{F}$-modules 
localement libres de rang fini ( fibr\'es vectoriels ) $Vect(F)$.\\

Un morphisme de champs $f : F \longrightarrow F'$ est projectif 
s'il existe une factorisation de $f$
$$\xymatrix{
F \ar[r]^{j} & \bf P\mit(V) \ar[r]^{p} & F' }$$
o\`u $j$ est une immersion ferm\'ee, et $p$ la projection d'un fibr\'e 
projectif associ\'e \`a un fibr\'e vectoriel $V$ sur $F'$.\\

\end{subsection}

\end{section}
\begin{section}{$G$-th\'eorie des Champs de Deligne-Mumford}
\hspace{5mm}

Dans cette section $k$ est un corps, et on suppose que $k$ contient 
les racines de l'unit\'e. Les champs alg\'ebriques sont des 
champs alg\'ebriques sur $(Sch/k)_{et}$.\\

Si $F$ est un champ alg\'ebrique, on dispose des cat\'egories
$Coh(F)$ et $Vect(F)$. Ce sont des cat\'egories exactes au 
sens de Quillen ( \cite[¤$2$]{q} ). On peut donc leur associer les spectres
de $K$-th\'eorie.

\begin{df}\label{d2a}
Les spectres de $K$-th\'eorie et de $G$-th\'eorie d'un champ alg\'ebrique $F$
sont d\'efinis par
$$K(F):=K(Vect(F))$$
$$G(F):=K(Coh(F))$$
On pose alors
$$K_{m}(F):=\pi_{m}(K(F))$$
$$G_{m}(F):=\pi_{m}(G(F))$$
\end{df}

La correspondance $F \mapsto K(F)$ ( resp. $F \mapsto G(F)$ ) est
un foncteur de la $2$-cat\'egorie des champs ( resp. champs et 
morphismes plats ), vers celles des spectres, morphismes de spectres et homotopies entre 
morphismes. Le produit tensoriel d\'efinit une structure de 
$K(F)$-module sur $G(F)$. 
Par restriction, on obtient deux pr\'efaisceaux en spectres 
$K$ et $G$ sur $F_{et}$.

\begin{df}\label{d2b}
Les spectres de $K$-th\'eorie et de $G$-th\'eorie \'etale sont d\'efinis par
$$K_{et}(F):=H(F_{et},K)$$
$$G_{et}(F):=H(F_{et},G)$$
On pose alors
$$K_{m,et}(F):=\pi_{m}(K_{et}(F))$$
$$G_{m,et}(F):=\pi_{m}(G_{et}(F))$$
\end{df}

Soit $p : F \rightarrow F'$ un morphisme propre de dimension 
cohomologique finie. En utilisant les arguments de \cite[$3.16.1$]{th}, on peut 
construire un morphisme d'images directes
$$p_{*} : G(F) \rightarrow G(F')$$
qui est  strictement fonctoriel pour la composition.\\
\hspace*{5mm} De la m\^eme fa\c{c}on, on dipose d'images r\'eciproques
$$f^{*} : G(F') \rightarrow G(F)$$
pour les morphismes $f : F \rightarrow F'$
de $Tor$-dimension finie.\\
\hspace*{5mm} On a \'egalement que $G_{et}(F)$ est fonctoriel covariant pour les 
images directes de morphismes propres repr\'esentables, et contravariant
pour les morphismes de $Tor$-dimension finie. La 
covariance pour les morphismes propres non n\'ecessairement 
repr\'esentables sera \'etudi\'ee plus loin.
Ces images directes et r\'eciproques sont li\'ees par les formules 
habituelles de transfert et de projection (\cite[¤$7$ $2.11$, $2.10$]{q} ).\\

\begin{subsection}{Propri\'et\'es g\'en\'erales}
\hspace{5mm}

La proposition suivante montre que les propri\'et\'es g\'en\'erales de la 
$G$-th\'eorie des champs sont les m\^emes que celles des sch\'emas.

\begin{prop}\label{p2a}
\begin{enumerate}
\item ( invariance topologique )
Soit $F$ un champ alg\'ebrique, et 
$j : F_{red} \hookrightarrow F$ l'immersion canonique de son sous-champ 
r\'eduit, alors
$$j_{*} : G(F_{red}) \rightarrow G(F)$$
est une \'equivalence faible.
\item ( localisation )
Si $j : F' \hookrightarrow F$ est un sous-champ ferm\'e de 
compl\'ementaire $i : F-F' \hookrightarrow F$, alors il existe une suite
exacte fonctorielle
$$\xymatrix{
G(F') \ar[r]^{j_{*}} & G(F) \ar[r]^{i^{*}} & G(F-F') }$$
\item ( axiome du fibr\'e projectif )
Si $p : V \rightarrow F$ est un fibr\'e vectoriel de rang $r+1$, $\pi : P \rightarrow F$ 
le fibr\'e projectif associ\'e, et $x:=\cal O\mit_{P}(1) \in K(P)$ le 
fibr\'e inversible canonique, alors le morphisme
$$\begin{array}{ccc}
\bigvee_{i=0}^{i=r}G(F) & \rightarrow & G(P) \\
\vee (a_{i}) & \mapsto & \sum_{i=0}^{i=r}\pi^{*}(a_{i}).x^{i}
\end{array}$$
est une \'equivalence faible.
\item ( homotopie )
Soit $p : V \rightarrow F$ un fibr\'e vectoriel. Alors le 
morphisme 
$$p^{*} : G(F) \rightarrow G(V)$$
est une \'equivalence faible.
\end{enumerate}
\end{prop}

\underline{\bf Preuve:} \rm Les preuves sont les m\^emes que pour le cas 
des sch\'emas. $\Box$\\ 

En localisant, on dispose de la m\^eme proposition pour le cas \'etale.

\begin{prop}\label{p2b}
\begin{enumerate}
\item ( invariance topologique ) Soit $F$ un champ alg\'ebrique, et
$j : F_{red} \hookrightarrow F$ l'immersion 
canonique de son sous-champ r\'eduit, alors
$$j_{*} : G_{et}(F_{red}) \rightarrow G_{et}(F)$$
est une \'equivalence faible.
\item ( localisation ) Si $j : F' \hookrightarrow F$ est un sous-champ ferm\'e de 
compl\'ementaire $i : F-F' \hookrightarrow F$, alors il existe une suite
exacte fonctorielle
$$\xymatrix{
G_{et}(F') \ar[r]^{j_{*}} & G_{et}(F) \ar[r]^{i^{*}} & G_{et}(F-F') }$$
\item ( axiome du fibr\'e projectif ) 
Si $p : V \rightarrow F$ est un fibr\'e vectoriel de rang $r+1$, $\pi : P \rightarrow F$ 
le fibr\'e projectif associ\'e, et $x:=\cal O\mit_{P}(1) \in K_{et}(P)$ le 
fibr\'e inversible canonique, alors le morphisme
$$\begin{array}{ccc}
\bigvee_{i=0}^{i=r}G_{et}(F) & \rightarrow & G_{et}(P) \\
\vee (a_{i}) & \mapsto & \sum_{i=0}^{i=r}\pi^{*}(a_{i}).x^{i}
\end{array}$$
est une \'equivalence faible.
\item ( homotopie ) Soit $p : V \rightarrow F$ un fibr\'e vectoriel. Alors le 
morphisme 
$$p^{*} : G_{et}(F) \rightarrow G_{et}(V)$$
est une \'equivalence faible.
\end{enumerate}
\end{prop}

Pour finir ces g\'en\'eralit\'es, rappelons le d\'evissage suivant la 
codimension du support \cite[$7.7$]{g2}.

\begin{prop}\label{p2c}
Soit $F$ un champ alg\'ebrique et $Coh(F)^{p}$, la cat\'egorie des 
faisceaux coh\'erents sur $F$ dont le support est de codimension 
sup\'erieur ou \'egal \`a $p$. Notons $F^{(p)}$ l'ensemble des points de 
codimension $p$ dans $F$, et $G(F)^{p}:=K(Coh(F)^{p})$. Alors les 
morphismes naturels d\'efinissent une suite exacte
$$\xymatrix{
G(F)^{p+1} \ar[r] & G(F)^{p} \ar[r] & \bigvee_{x \in 
F^{(p)}}G(\widetilde{x}) }$$ 
\end{prop}

\end{subsection}
 
\begin{subsection}{Th\'eor\`emes de Descente}
\hspace{5mm}

Nous venons de voir que la $G$-th\'eorie des champs alg\'ebriques 
poss\`edait beaucoup de propri\'et\'es analogues au cas des sch\'emas. 
Cependant le th\'eor\`eme de descente \'etale de la $G$-th\'eorie rationnelle 
s'av\`ere faux en g\'en\'eral. Nous allons voir qu'une g\'en\'eralisation plus 
faible existe tout de m\^eme. Elle nous permettra de nous ramener 
souvent au cas des champs quotients par des groupes finis.
Nous d\'emontrerons par la suite un th\'eor\`eme de descente homologique. Ce r\'esultat 
nous permettra de d\'efinir les images directes pour la $G$-th\'eorie \'etale.\\ 

Nous adopterons la notation suivante : si $E$ est un pr\'efaisceau en 
spectres on notera $E_{\bf Q}:=E \otimes \bf Q$. C'est le pr\'efaisceau 
localis\'e de $E$ suivant les \'equivalences faibles rationnelles.\\

\begin{subsubsection}{Descente Etale}

\begin{thm}\label{t2a}
Soit $p : F \rightarrow M$ la projection d'un champ alg\'ebrique sur son
espace de modules. Alors $p_{*}G_{\bf Q}$ est flasque sur $M_{et}$.
\end{thm}

Remarquons que si $F$ est un sch\'ema, le th\'eor\`eme 
redonne le th\'eor\`eme de descente \'etale de la $G$-th\'eorie rationnelle ( 
\cite[$11.10$]{th} ). \\

\underline{\bf Preuve:} \rm En utilisant le "lemme des cinq", on sait 
que s'il existe une suite exacte
$$\xymatrix{ E \ar[r] & F \ar[r] & G }$$
sur $M_{et}$, $F$ est flasque si $G$ et $E$ le sont. \\
\hspace*{5mm} Utilisons alors la proposition \ref{p2c}. On dispose d'une suite 
exacte sur $M_{et}$
$$\xymatrix{
p_{*}G_{\bf Q}^{p+1} \ar[r] & p_{*}G_{\bf Q}^{p} \ar[r] & \bigvee_{x \in 
F^{(p)}}(i_{x})_{*}(p_{x})_{*}G_{\bf Q} }$$
o\`u $i_{x} : Spec k(x) \rightarrow M$
est le morphisme associ\'e au point $x$, et $p_{x} : 
\widetilde{x} \rightarrow Spec k(x)$ la projection naturelle. En raisonnant par 
r\'ecurrence descendante sur $p$, il nous suffit de montrer que le terme 
de droite est flasque sur $M_{et}$, et donc que pour chaque point $x$ 
de $M$, le pr\'efaisceau $(i_{x})_{*}(p_{x})_{*}G_{\bf Q}$ est flasque sur 
$M$. Comme les images directes pr\'eservent les pr\'efaisceaux flasques, 
il suffit de montrer que $(p_{x})_{*}G_{\bf Q}$ est fasque sur $(Spec 
k(x))_{et}$. Ainsi on se ram\`ene \`a d\'emontrer le th\'eor\`eme dans le cas o\`u $M=Spec K$ est le 
spectre d'un corps. Comme $M_{et}$ est alors le site galoisien de $K$, 
le th\'eor\`eme provient du lemme suivant.

\begin{lem}\label{l2a}
Soit $F \rightarrow Spec K$ une gerbe sur un corps $K$, $K^{sp}$ une 
cl\^oture s\'eparable de $K$, et $H=Gal(K^{sp}/K)$. On note $F^{sp} 
\rightarrow Spec K^{sp}$ la gerbe induite. Alors le morphisme canonique 
$q : F^{sp} \rightarrow F$, induit une \'equivalence faible
$$q^{*} : G_{\bf Q}(F) \rightarrow G_{\bf Q}(F^{sp})^{H}$$
\end{lem}

\underline{\bf Preuve:} \rm 
On commence par \'eclaircir les notations. 
Le groupe $H$ op\`ere par automorphismes sur $Spec K^{sp}$, et par 
fonctorialit\'e, on a une op\'eration $H \rightarrow 
Aut_{F}(F^{sp})$. Donc $H$ op\`ere sur l'espace $G(F^{sp})_{\bf Q}$, que l'on voit alors comme un 
$H$-diagramme. Dans ce cas 
$$G(F^{sp})_{\bf Q}^{H}:=Holim_{H}G(F^{sp})_{\bf Q}$$
est l'espace des invariants homotopiques de $H$.\\
\hspace*{5mm} Comme $H$ est profini, et que $G(F^{sp})_{\bf Q}$ est un 
spectre rationnel, la suite spectrale de
 Bousfield-Kan ( \cite[Ch. $XI$ \S $7$]{bk} ) d\'eg\'en\`ere et donne des 
isomorphismes canoniques $\pi_{m}(G(F^{sp})_{\bf Q}^{H}) \simeq 
G_{m}(F^{sp})_{\bf Q}^{H}$. De plus par continuit\'e du foncteur $G$, on a 
$$G_{m}(F^{sp})_{\bf Q}^{H} \simeq Colim_{K 
\hookrightarrow K_{i} \; galois \; fini}G_{m}(F_{i})_{\bf Q}^{H_{i}}$$
avec $H_{i}=Gal(K_{i}/K)$, et $F_{i}=F\times_{Spec K}Spec 
K_{i}$.
Il reste \`a montrer que le morphisme naturel $G(F)_{\bf Q} \rightarrow 
G(F_{i})_{\bf Q}^{H_{i}}$ est une \'equivalence faible pour tout $i$. \\
\hspace*{5mm} Soit $p : F_{i} \rightarrow F$ la projection, et $m_{i}$ l'ordre de $H_{i}$.
Alors la formule de projection implique que 
$$p_{*} \circ p^{*} = \times m_{i} : G(F)_{\bf Q} \rightarrow G(F)_{\bf Q}$$
 Et de m\^eme 
$$p^{*} \circ p_{*} = \sum_{h \in H_{i}} h^{*} \sim \times m_{i} : 
G(F_{i})_{\bf Q}^{H_{i}} \rightarrow G(F_{i})_{\bf Q}^{H_{i}}$$
Ainsi, $\frac{1}{m_{i}}p_{*}$ est un inverse homotopique de $p^{*}$. $\Box$\\

\begin{cor}\label{c2a}
Soit $F$ un champ alg\'ebrique, et $p : F \rightarrow M$ son espace de modules. 
\begin{enumerate}
\item Pour tout recouvrement \'etale $U \rightarrow M$, le morphisme naturel
$$G(F)_{\bf Q} \rightarrow \Check{H}(U/M,p_{*}G_{\bf Q})$$
est une \'equivalence faible.
\item Soit $X \rightarrow M$ un rev\^etement galoisien de groupe $H$, 
et $p : F_{X}:=F\times_{M}X \rightarrow F$ le rev\^etement induit. Alors les 
morphismes naturels
$$p^{*} : G(F)_{\bf Q} \rightarrow G(F_{X})^{H}_{\bf Q}$$
$$p_{*} : G(F_{X})_{\bf Q}^{H} \rightarrow G(F)_{\bf Q}$$
sont des \'equivalences faibles.
\item Si $F$ est lisse, alors il existe une \'equivalence faible 
naturelle
$$H(M_{et},p_{*}K_{\bf Q}) \simeq G(F)_{\bf Q}$$
\end{enumerate}
\end{cor}

\underline{\bf Preuve:} \rm Les points $(1)$ et $(2)$ proviennent du 
th\'eor\`eme de descente ( \ref{desc} ) appliqu\'e au pr\'efaisceau flasque $p_{*}G_{\bf Q}$ 
sur $M_{et}$, et de la formule de projection $p_{*}\circ p^{*}=\times 
o(H)$ pour le $(2)$. \\
\hspace*{5mm} Pour d\'emontrer $(3)$, il suffit d'apr\`es le th\'eor\`eme de d\'emontrer que
$$p_{*}K_{\bf Q} \rightarrow p_{*}G_{\bf Q}$$
est une \'equivalence faible de pr\'efaisceaux en spectres sur $M_{et}$. 
En localisant sur $M_{et}$, on se ram\`ene au cas o\`u $F=[X/H]$ est un 
champ quotient d'un sch\'ema lisse $X$, par un groupe fini $H$. Mais 
alors on sait que tout faisceau coh\'erent sur $F$ admet une r\'esolution 
par des fibr\'es vectoriels ( \cite[$5.3$]{th2} ), et donc que le morphisme naturel
$$K(F) \rightarrow G(F)$$ est une \'equivalence faible.$\Box$\\

\end{subsubsection}

\begin{subsubsection}{Descente Homologique}
\hspace{5mm}

Nous allons maintenant d\'emontrer un r\'esultat de descente de la 
$G$-th\'eorie \'etale. Il va nous permettre d'\'etendre la d\'efinition des images 
directes aux morphismes propres quelconques. La m\'ethode de 
construction est tir\'ee de \cite{g3}.\\

Soit $F$ un champ alg\'ebrique et $p : X \rightarrow F$ un morphisme 
propre surjectif, avec $X$ un sch\'ema. Le nerf de $X$ sur $F$ est un 
sch\'ema simplicial \`a faces propres $\cal N\mit(X/F)$. On peut donc en prendre l'image 
par le foncteur covariant $G$, et obtenir ainsi un spectre 
simplicial not\'e $G(\cal N\mit(X/F))$. La colimite homotopique de ce 
spectre simplicial sera not\'ee
$$G(X/F):=Hocolim_{\Delta^{op}}G(\cal N\mit(X/F))$$

\begin{thm}\label{t2b}
Soit $X$ un sch\'ema et $p : X \rightarrow F$ un morphisme propre 
surjectif. Alors le morphisme naturel
$$p_{*} : G(X/F)_{\bf Q} \rightarrow G_{et}(F)_{\bf Q}$$
est une \'equivalence faible.
\end{thm}

\underline{\bf Preuve:} \rm En raisonnant par r\'ecurrence sur la 
dimension de $F$ et \`a l'aide de \ref{p2a}, on peut se restreindre \`a
d\'emontrer le th\'eor\`eme pour un ouvert de Zariski de $F$. Ainsi on se 
ram\`ene au cas o\`u $F$ est une gerbe sur $M$, telle qu'il existe un 
rev\^etement galoisien $M' \rightarrow M$ fini et une 
\'equivalence
$$F':=F\times_{M}M' \simeq [V/H]$$ avec $V$ un sch\'ema et $H$ un groupe 
fini op\'erant sur $V$.  Notons $X':=X\times_{F}F'$.  \\
\hspace*{5mm} Rappelons que si $E$ est un spectre muni d'une action d'un groupe 
$H$, on note $E_{H}:=Hocolim_{H}E$. \\
\hspace*{5mm} Consid\'erons alors le diagramme homotopiquement commutatif d'images directes suivant 
$$\xymatrix{ G(X/F)_{\bf Q} \ar[r] & G_{et}(F)_{\bf Q} \\
(G(X'/F')_{\bf Q})_{H} \ar[u]Ê\ar[r] & (G_{et}(F')_{\bf Q})_{H} \ar[u] }$$
Comme les fl\`eches verticales sont des \'equivalences faibles ( d'apr\`es
\ref{c2a} ), on se ram\`ene \`a d\'emontrer que $(G(X'/F')_{\bf Q})_{H} \rightarrow 
(G_{et}(F')_{\bf Q})_{H}$ est une \'equivalence faible.  Mais comme les 
co-invariants homotopiques pr\'eservent les \'equivalences faibles, il 
nous suffit de d\'emontrer le th\'eor\`eme pour $X' \rightarrow F'$.  En 
utilisant le m\^eme argument avec le rev\^etement galoisien $V 
\rightarrow F'$, il suffit de d\'emontrer le th\'eor\`eme dans le cas o\`u $F$ est un 
sch\'ema.\\
\hspace*{5mm} Comme pr\'ec\'edemment on peut se restreindre \`a un ouvert g\'en\'erique de 
$F$.  On peut alors supposer que $X \rightarrow F$ poss\`ede une section 
apr\`es un changement de base fini $F' \rightarrow F$.  On peut m\^eme supposer que
$F' \rightarrow F$ se factorise en
$$\xymatrix{ F' \ar[r]^{a} & F'' \ar[r]^{b} & F }$$
avec $a$ un morphisme radiciel, et $b$ un rev\^etement 
galoisien de groupe $H$. Notons $c : X':=X\times_{F}F' \rightarrow F'$ et 
$d : X'':=X\times_{F}F'' \rightarrow F''$.
On dispose alors du diagramme homotopiquement commutatif suivant
$$\xymatrix{
G(X/F)_{\bf Q} \ar[r]^{p_{*}} & G(F)_{\bf Q} \\
G(X''/F'')_{\bf Q} \ar[u] \ar[r]^{d_{*}} & G(F'')_{\bf Q} 
\ar[u]^{b_{*}} \\
G(X'/F')_{\bf Q} \ar[u] \ar[r]^{c_{*}} & G(F')_{\bf Q} 
\ar[u]^{a_{*}} }$$
Comme $a$ est radiciel, on sait que $a_{*} : G(F')_{\bf Q} \rightarrow 
G(F'')_{\bf Q}$ est une \'equivalence faible. Il en est de m\^eme du 
morphisme induit $G(X'/F')_{\bf Q} \rightarrow G(X''/F'')_{\bf Q}$. 
De plus, comme le morphisme $X' \rightarrow F'$ poss\`ede une section, 
$c_{*} : G(X'/F')_{\bf Q} \rightarrow G(F')_{\bf Q}$ est une 
\'equivalence faible ( \cite[$4-1$ $(I)$]{g3} ). Ceci implique que le th\'eor\`eme est 
vrai pour $X'' \rightarrow F''$. En prenant les co-invariants 
homotopiques on obtient un nouveau diagramme qui commute \`a homotopie 
pr\`es
$$\xymatrix{
G(X/F)_{\bf Q} \ar[r]^{p_{*}} & G(F)_{\bf Q} \\
(G(X''/F'')_{\bf Q})_{H} \ar[u] \ar[r]^{d_{*}} & (G(F'')_{\bf Q})_{H} 
\ar[u]^{b_{*}} }$$
De nouveau les fl\`eches verticales sont des \'equivalences faibles. Le 
th\'eor\`eme \'etant vrai pour $X''/F''$, il est vrai pour $X/F$. $\Box$.\\

Passons \`a la construction des images directes en $G$-th\'eorie \'etale. \\
\hspace*{5mm} Soit $f : F \rightarrow F'$ un morphisme propre. A l'aide de 
\cite[$4.12$]{dm}, choisissons un 
sch\'ema $X$ et un morphisme propre surjectif
$$p : X \rightarrow F$$
On dispose alors d'un diagramme de morphismes propres
$$\xymatrix{
\cal N\mit (X/F) \ar[d]_{p} \ar[rd]^{q} & \\
F \ar[r]^{f} & F' }$$
En appliquant le foncteur covariant $(G_{et})_{\bf Q}$, on obtient un 
diagramme de spectres
$$\xymatrix{
G(X/F)_{\bf Q} \ar[d]_{p_{*}} \ar[rd]^{q_{*}} & \\
G_{et}(F)_{\bf Q}  & G_{et}(F')_{\bf Q} }$$
D'apr\`es le th\'eor\`eme \ref{t2b}, $p_{*}$ est une \'equivalence faible. On 
dispose donc d'un morphisme bien d\'efini \`a homotopie pr\`es
$$f_{*}=q_{*}\circ (p_{*}^{-1}) : G_{et}(F)_{\bf Q} \rightarrow 
G_{et}(F')_{\bf Q}$$
Soit $p' : X' \rightarrow F$ est un autre morphisme propre surjectif, 
et $Z=X\times_{F}X'$. Alors, les diagrammes suivants 
commutent \`a homotopie pr\`es
$$\xymatrix{
G(Z/F)_{\bf Q} \ar[r] \ar[d] & G(X/F)_{\bf Q} \ar [d]^{p_{*}} \\
G(X'/F)_{\bf Q} \ar[r]^{(p')_{*}} & G_{et}(F)_{\bf Q} }$$
$$\xymatrix{
G(Z/F)_{\bf Q} \ar[r] \ar[d] & G(X/F)_{\bf Q} \ar [d]^{q_{*}} \\
G(X'/F)_{\bf Q} \ar[r]^{(q')_{*}} & G_{et}(F')_{\bf Q} }$$
Ainsi, le morphisme $f_{*} : G_{et}(F)_{\bf Q} \rightarrow G_{et}(F')_{\bf Q}$ 
est ind\'ependant \`a homotopie pr\`es du choix de $p : X \rightarrow F$. 
Il est alors facile de voir que $F \rightarrow G_{et}(F)$ est un 
foncteur covariant \`a homotopie pr\`es pour tous les morphismes propres.

\begin{lem}\label{l2b}
Soit $X$ un sch\'ema et $H$ un groupe fini op\'erant sur $X$. Alors 
la projection $p : X \rightarrow X/H$ induit des \'equivalences faibles
$$p_{*} : (G(X)_{\bf Q})_{H} \rightarrow G(X/H)_{\bf Q}$$
$$p_{*} : G(X)_{\bf Q}^{H} \rightarrow G(X/H)_{\bf Q}$$
\end{lem}

\underline{\bf Peuve:} \rm Il suffit de consid\'erer le cas o\`u $X$ r\'eduit.\\

On consid\`ere le diagramme commutatif
$$\xymatrix{
G(X)_{\bf Q}^{H} \ar[r] \ar[rd]_{p_{*}} & (G(X)_{\bf Q})_{H} 
\ar[d]^{p_{*}} \\
 & G(X/H)_{\bf Q} }$$
Comme le groupe $H$ est fini, et que l'on est \`a coefficients 
rationnels, le morphisme horizontal est une \'equivalence faible. Il 
suffit donc de d\'emontrer la seconde assertion du lemme. Soit $j : Y 
\hookrightarrow X$ une immersion ferm\'e \'equivariante, tel que $i : X-Y 
\hookrightarrow X$ soit dense dans $X$, et que l'action de $H$ soit sans 
points fixes 
sur $X-Y$. On consid\`ere le morphisme de suites exactes
$$\xymatrix{
G(Y)_{\bf Q}^{H} \ar[r]^{j_{*}} \ar[d]_{p_{*}} & G(X)_{\bf Q}^{H} 
\ar[r]^{i^{*}} \ar[d]_{p_{*}} & G(X-Y)_{\bf Q}^{H} \ar[d]_{p_{*}} \\
G(Y/H)_{\bf Q} \ar[r]^{j_{*}} & G(X/H)_{\bf Q} \ar[r]^{i^{*}} & 
G((X-Y)/H)_{\bf Q} }$$
En raisonnant par r\'ecurrence sur la dimension de $X$, on peut se 
restreindre au cas o\`u l'action de $H$ est sans points fixes sur $X$. 
Le corollaire \ref{c2a} $(2)$ montre alors que 
$p_{*} : G(X)^{H}_{\bf Q} \rightarrow G(X/H)_{\bf Q}$ est une 
\'equivalence faible. $\Box$\\

\begin{cor}\label{c2b}
Soit $p : F \rightarrow M$ la projection d'un champ alg\'ebrique sur son espace de 
modules. Alors le morphisme
$$p_{*} : G_{et}(F)_{\bf Q} \rightarrow G(M)_{\bf Q}$$
est une \'equivalence faible.
\end{cor}

\underline{\bf Preuve} \rm On consid\`ere le morphisme de pr\'efaisceaux 
en spectres sur $M_{et}$
$$p_{*} : p_{*}(G_{et})_{\bf Q} \rightarrow G_{\bf Q}$$
On sait que le membre de droite est flasque sur $M_{et}$ ( \ref{t2a} 
pour un espace alg\'ebrique ). Le membre de gauche est  l'image
directe d'un pr\'efaisceau flasque, donc est encore flasque. On peut 
donc localiser sur $M_{et}$ et se ramener au cas o\`u $F=[X/H]$ est le 
champ quotient d'un sch\'ema par un groupe fini $H$. \\
\hspace*{5mm} Prenons le morphisme naturel $q : X \rightarrow F$, et
notons $r : X \rightarrow X/H=M$. On peut alors consid\`erer le diagramme 
homotopiquement commutatif suivant
$$\xymatrix{
G(X/F)_{\bf Q} \ar[r]^{r_{*}} \ar[d]_{q_{*}} & G(M)_{\bf Q} \\
G_{et}(F)_{\bf Q} \ar[ur]_{p_{*}} }$$
Par le th\'eor\`eme \ref{t2b}, $q_{*}$ est une \'equivalence faible.
Il nous suffit donc de montrer que $r_{*}$ est une \'equivalence faible. 
Or le nerf $\cal N\mit (X/F)$ est canoniquement isomorphe au 
classifiant $BHX$ de l'action de $H$ sur $X$. Il existe donc une 
\'equivalence faible naturelle
$$(G(X)_{\bf Q})_{H} \rightarrow G(X/F)_{\bf Q}$$
qui rend homotopiquement commutatif le diagramme
$$\xymatrix{
(G(X)_{\bf Q})_{H} \ar[r] \ar[d]_{r_{*}} & G(X/F)_{\bf Q} \ar[dl]_{r_{*}} \\
G_{et}(M)_{\bf Q} & }$$
Le lemme pr\'ec\'edent permet de conclure que $r_{*}$ est une 
\'equivalence faible. $\Box$\\

Il est \`a noter que ces morphismes d'images directes ne sont pas 
d\'efinis \`a coefficients entiers, tout au moins lorsque les morphismes ne 
sont pas repr\'esentables. Ils sont l'analogue des images directes de 
cycles alg\'ebriques d\'efinis dans \cite{vi2}. Ainsi le corollaire 
pr\'ec\'edent est un analogue au fait que les groupes de Chow rationnels d'un champ 
alg\'ebrique sont isomorphes \`a ceux de son espace de modules. 

\end{subsubsection}

\begin{subsubsection}{Enveloppes de Chow}
\hspace{5mm}

Dans le cas des sch\'emas, on sait qu'une enveloppe de Chow $p : Z 
\rightarrow X$ permet de retrouver la $G$-th\'eorie de $X$ ( 
\cite{g3} ). Si maintenant $F$ est un champ alg\'ebrique, on verra
qu'il existe une notion d'enveloppe de Chow rationnelle permettant de calculer la
$G$-th\'eorie rationnelle de $F$, en fonction de la $G$-th\'eorie de 
certaines gerbes triviales. \\

\begin{df}\label{d2c}
Soit $F$ un champ alg\'ebrique. Un morphisme propre et repr\'esentable de champs
$$p : F_{0} \rightarrow F$$
est une enveloppe de Chow rationnelle, si 
pour tout point $x$ de $F$, le morphisme induit sur la gerbe r\'esiduelle
$$p^{-1}(\widetilde{x}) \longrightarrow \widetilde{x}$$
poss\`ede une section apr\`es une extension finie de $k(x)$.
\end{df}  

Remarquons par exemple, que le morphisme naturel $X \rightarrow 
[X/H]$ est une enveloppe de Chow rationnelle, si et seulement si $H$ op\`ere 
sans points fixes sur $X$.\\

La condition de la d\'efinition est \'equivalente \`a dire que pour tout 
sous-champ ferm\'e int\`egre $F'$ de $F$, il existe un sous-champ ferm\'e int\`egre
$F_{0}'$ de $F_{0}$ au-dessus de $F'$, tel que la restriction 
$$p : F'_{0} \rightarrow F'$$
admette g\'en\'eriquement une section apr\`es un changement de base fini de 
l'espace de modules de $F'$.\\

Si $F' \rightarrow F$ est un morphisme propre de champs, on dispose de son 
nerf $\cal N\mit(F'/F)$, qui est un champ simplicial \`a faces propres. En 
appliquant le foncteur covariant $G$, on en d\'eduit un spectre simplicial 
$G(\cal N\mit(F'/F))$. La colimite homotopique de ce diagramme sera 
not\'ee $G(F'/F)$.\\ 

\begin{prop}\label{p2d}
Si $p : F_{0} \rightarrow F$ est une enveloppe de Chow rationnelle, alors le morphisme 
naturel
$$p_{*} : G(F_{0}/F)_{\bf Q} \rightarrow G(F)_{\bf Q}$$
est une \'equivalence faible. 
\end{prop}

\underline{\bf Preuve:} \rm Comme dans la preuve de \ref{t2b}, on peut se 
restreindre au cas o\`u $F_{0} \longrightarrow F$ poss\`ede une section apr\`es un
changement de base de l'espace de modules de $F$ par un rev\^etement 
galoisien. En utilisant \ref{c2a} $(2)$ on se ram\`ene au cas o\`u $p$ poss\`ede 
une section $s : F \rightarrow F_{0}$. On sait alors que 
$s_{*} : G(F)_{\bf Q} \rightarrow G(F_{0}/F)_{\bf Q}$ est un inverse homotopique 
de $p_{*}$ ( \cite[$4-1$ $(I)$]{g3} ).
$\Box$\\

Remarquons que les arguments de la preuve de  
\cite[$4-1$]{g3}, entrainent que la proposition se g\'en\'eralise au cas des 
hyper-enveloppes de Chow rationnelles. 

\end{subsubsection}

\end{subsection}

\begin{subsection}{Th\'eor\`eme de D\'evissage}
\hspace{5mm}

Dans cette section  nous allons montrer que les espaces de 
$G$-th\'eorie, se d\'ecrivent \`a l'aide de la 
$G$-th\'eorie \'etale du champ des ramifications. Pour cela, nous allons 
g\'en\'eraliser certains r\'esultats de Vistoli ( \cite{vi} ) au cas des 
champs alg\'ebriques.\\

Nous aurons besoin d'une hypoth\`ese suppl\'ementaire. Elle assure par 
exemple que les morphismes propres sont tous de dimension 
cohomologique finie.

\begin{hyp}\label{h}
Si $F$ est un champ alg\'ebrique, on supposera
que les ordres d'inertie de tous les points de $F$ sont premiers \`a la 
caract\'eristique de $k$.
\end{hyp}

Pour ce chapitre, on notera $\Lambda=\bf Q\mit(\mu_{\infty})$, avec 
$\mu_{\infty}$ le sous-groupe de $k^{*}$ des racines de l'unit\'e, que l'on
identifiera \`a un sous-groupe de $\bf C\mit^{*}$, par un plongement fix\'e
$$\mu_{\infty} \hookrightarrow \bf C\mit^{*}$$ 
Si $E$ est un spectre, $E_{\Lambda}$ d\'esignera $E\otimes_{\bf Z}\bf 
Q\mit(\mu_{\infty})$.\\

\begin{thm}\label{t2d}
Soit $F$ un champ alg\'ebrique lisse. Alors il existe une \'equivalence 
faible de spectres en $\Lambda$-alg\`ebres
$$\phi_{F} : G(F)_{\Lambda} \stackrel{\simeq}{\longrightarrow} 
G_{et}(I_{F})_{\Lambda}$$
De plus, $\phi$ commute avec les images r\'eciproques, et les produits.
\end{thm}

\underline{\bf Preuve:} \rm Commen\c{c}ons par d\'efinir un morphisme
$$\rho : K(I_{F}) \rightarrow K(I_{F})_{\Lambda}$$
Soit $V$ un fibr\'e vectoriel sur $I_{F}$. Fixons nous $\zeta$ 
une racine de l'unit\'e dans $k$.\\

Soit $f : X \rightarrow I_{F}$ une section au-dessus d'un sch\'ema $X$. 
Elle est d\'etermin\'ee par une section $s : X \rightarrow F$, et 
un automorphisme $a \in \pi_{1}(F(X),s)$. Le fibr\'e $f^{*}V$ sur 
$X$ est alors muni d'une action de $<a>$. D'apr\`es \ref{h}, cette 
action est diagonalisable. On obtient donc une d\'ecomposition canonique
$$f^{*}V \simeq f^{*}V^{(\zeta)}\oplus W$$
o\`u $f^{*}V^{(\zeta)}$ est le sous-fibr\'e vectoriel de $f^{*}V$ sur lequel 
$a$ op\`ere par multiplication par $\zeta$.\\
\hspace*{5mm} Donnons-nous maintenant deux sections
$$f : X \rightarrow I_{F}$$
$$g : Y \rightarrow I_{F}$$ 
donn\'ees respectivement par $(s,a)$ et $(t,b)$, un morphisme $p : Y 
\rightarrow X$, et une homotopie $h : f \circ p \Rightarrow g$. 
Le cocyle de $V$
$$\theta_{f,g,h} : (p\circ f)^{*}V \simeq g^{*}V$$
\'etant compatible avec les actions de $<a>$ et $<b>$, on 
obtient ainsi un cocycle induit 
$$\theta_{f,g,h} : (p\circ f)^{*}V^{(\zeta)} \simeq g^{*}V^{(\zeta)}$$ 
De cette fa\c{c}on on d\'efinit un fibr\'e vectoriel $V^{(\zeta)}$ sur $I_{F}$. De 
plus, Il est clair que les foncteurs
$V \mapsto V^{(\zeta)}$
sont exacts. On pose alors
$$\begin{array}{cccc}
\rho : & K(I_{F}) & \longrightarrow & K(I_{F})_{\Lambda} \\
       &    V     & \mapsto         & \sum_{\zeta \in \mu_{\infty}} 
       \zeta.V^{(\zeta)} 
\end{array}$$
C'est un morphisme de spectres en anneaux, qui est de plus fonctoriel pour 
les images r\'eciproques.
Soit $\pi : I_{F} \longrightarrow F$ la projection, et $can : K(I_{F}) 
\longrightarrow K_{et}(I_{F})$ le morphisme canonique. On consid\`ere 
$$can \circ \rho \circ \pi^{*} : K(F)_{\Lambda} \rightarrow 
K_{et}(I_{F})_{\Lambda}$$
En localisant sur $M_{et}$, on obtient 
un morphisme de pr\'efaisceaux en spectres sur $M_{et}$
$$\phi_{F} : p_{*}K_{\Lambda} \longrightarrow 
p_{*}\pi_{*}(K_{et})_{\Lambda}$$
Et par le corollaire \ref{c2a} $(3)$, le morphisme cherch\'e 
$$\phi_{F} : G_{*}(F)_{\Lambda} \longrightarrow G_{*,et}(I_{F})_{\Lambda}$$ 
Par construction, il est compatible aux produits et aux images 
r\'eciproques. Enfin, il est clair que lorsque $F=[X/H]$, on obtient le morphisme
$\phi$ d\'efini dans \cite[$4.1$]{vi}. Ainsi $\phi_{F}$ est une 
\'equivalence faible. $\Box$\\

Remarquons que m\^eme si $F$ n'est pas lisse, il existe toujours un morphisme
$$\phi_{F} : K(F)_{\Lambda} \longrightarrow K_{et}(I_{F})_{\Lambda}$$
mais qui n'est pas, \`a priori, une \'equivalence faible. \\

Enfin, si $F$ est une gerbe ( \'eventuellement singuli\`ere ), il existe encore 
une \'equivalence faible canonique de spectres en $\Lambda$-modules
$$\phi_{F} : G(F)_{\Lambda} \longrightarrow G_{et}(I_{F})_{\Lambda}$$
En effet, comme $\pi : I_{F} \longrightarrow F$ est \'etale, on peut 
appliquer la construction de la preuve pr\'ec\'edente. 

\end{subsection}

\end{section}
\begin{section}{Le Th\'eor\`eme de Grothendieck-Riemann-Roch}
\hspace{5mm}

On continue \`a supposer que $k$ contient les racines de l'unit\'e, et 
que les champs v\'erifient \ref{h}.

\begin{subsection}{Cohomologie des Champs Alg\'ebriques}
\hspace{5mm}

On se place dans $(SchQP/k)_{et}$, 
le site des sch\'emas quasi-projectifs sur $k$, muni de la topologie \'etale. 
Par la suite tous les sch\'emas sont quasi-projectifs. \\

Donnons-nous une th\'eorie cohomologique avec images 
directes. C'est \`a dire les donn\'ees suivantes.

\begin{itemize}
\item pour tout entier $i$, un pr\'efaisceau en spectres ab\'eliens $\cal H\mit^{i}$ 
( i.e. associ\'e par la construction de Dold-Puppe \`a un complexe de 
pr\'efaisceaux ab\'eliens )
flasque sur le site pr\'ec\'edent, et muni d'un produit 
$$\cal H\mit^{i} \wedge \cal H\mit^{j} \longrightarrow \cal H\mit^{i+j}$$ 
homotopiquement associatif et commutatif. On notera 
$$\cal H\mit:=\bigvee_{i \in \bf Z}\cal H\mit^{i}$$ 
le pr\'efaisceau en spectres en anneaux gradu\'es associ\'e.
\item pour chaque entier $i$, et chaque 
sch\'ema $X$, des spectres ab\'eliens $\cal H\mit'_{i}(X)$. Pour chaque 
morphisme propre $f : X \rightarrow Y$ de sch\'emas irr\'eductibles 
un morphisme de spectres ab\'eliens
$$f_{*} : \cal H\mit'_{i}(X) \longrightarrow \cal H\mit'_{i+dp}(Y)$$
avec $p$ la dimension relative de $X$ sur $Y$, et
$d$ un entier fix\'e, \'egal \`a $1$ ou $2$. On demande que ces 
images directes s'\'etendent en un foncteur 
$$X \mapsto \cal H\mit'(X):=\bigvee_{i \in \bf Z} \cal H\mit'_{i}(X)$$
covariant pour les morphismes propres.
\item Pour chaque sch\'ema $X$, une structure de $\cal H\mit(X)$-module 
gradu\'e sur $\cal H\mit'(X)$. De plus, 
si $f : X \longrightarrow Y$ est un morphisme propre, et $x \in 
\cal H\mit'(X)$, $y \in \cal H\mit(Y)$, alors on a 
$$f_{*}(f^{*}(y).x)=y.f_{*}(x)$$
\item pour chaque morphisme \'etale $f : X \rightarrow Y$, un morphisme 
de spectres ab\'eliens
$$f^{*} : \cal H\mit'_{i}(Y) \longrightarrow \cal H\mit'_{i}(X)$$
qui fait de $\cal H\mit'$ un $\cal H\mit$-module
pour les morphismes \'etales.
\item si 
$$\xymatrix{
X' \ar[r]^{u} \ar[d]_{q} & X \ar[d]^{p} \\
X \ar[r]_{v} & Y }$$
est cart\'esien
avec $p$ propre et $u$ \'etale, alors 
$$v^{*} \circ p_{*} = q_{*} \circ u^{*}$$
\item si $j : Y \hookrightarrow X$ est une immersion ferm\'ee, et 
$i : U \hookrightarrow X$ l'immersion ouverte compl\'ementaire, alors 
la suite suivante
$$\xymatrix{
\cal H\mit'_{i}(Y) \ar[r]^{j_{*}} & \cal H\mit'_{i+dp}(X) 
\ar[r]^{i^{*}} & \cal H\mit'_{i+dp}(U) }$$
est canoniquement une suite exacte.
\item si $X$ est un sch\'ema lisse, alors il existe une
\'equivalence faible
$$p_{X} : \cal H\mit^{i}(X)  \longrightarrow  \cal H\mit'_{i}(X)$$
compatible avec les produits et les images r\'eciproques.
\item Il existe un morphisme de pr\'efaisceaux  
simpliciaux 
$$C_{1} : B\bf G\mit_{m}  \longrightarrow \cal H\mit^{d}_{[0]}$$ 
Et donc un morphisme flexible ( \cite{s} ) sur les champs associ\'es
$$C_{1} : Pic \longrightarrow \cal H\mit^{d}_{[0]}$$
\item soit $\pi : P \longrightarrow X$ le fibr\'e projectif associ\'e \`a
un fibr\'e vectoriel de rang $r+1$, et $x=C_{1}(\cal O\mit_{P}(1))$. 
Alors le morphisme 
$$\begin{array}{ccc}
\bigvee_{i=0}^{i=r} \cal H\mit(X) & \longrightarrow & \cal H\mit(P) \\
\vee a_{i} & \mapsto & \sum_{i}\pi^{*}(a_{i}).x 
\end{array}$$
est une \'equivalence faible.
\end{itemize}

Avant d'aller plus loin, donnons quelques exemples. 

\begin{itemize}

\item \underline{La th\'eorie de Gersten:} Notons $\cal K\mit_{i}$ le 
pr\'efaisceau sur $(SchQP/k)_{et}$
$$X \mapsto K_{i}(X)$$
On consid\`ere les pr\'efaisceaux en spectres $K(\cal K\mit_{i},i)\otimes 
\bf Q$, et
on prend pour $\cal H\mit^{i}$ les pr\'efaisceaux fibrants associ\'es.\\
\hspace*{5mm} Soit 
$$\cal R\mit^{i} : \bigoplus_{x \in X^{(0)}}K_{i}(k(x)) \longrightarrow \bigoplus_{x \in 
X^{(1)}}K_{i-1}(k(x)) \longrightarrow \dots \bigoplus_{x \in X^{(i)}}K_{0}(k(x))$$
le complexe de Gersten de codimension $i$ de $X$, concentr\'e en 
degr\'e $[-i,0]$. On pose 
$$\cal H\mit_{i}'(X):=K(\cal R\mit^{i},0)\otimes \bf Q$$
Alors $\cal H$ et $\cal H\mit'$ v\'erifient les axiomes ci-dessus ( 
\cite[$1.4$ $(ii)$]{g} ), avec $d=1$.

\item \underline{La cohomologie de De Rham:} Supposons $k$ de 
caract\'eristique nulle. \\
\hspace*{5mm} Pour un sch\'ema $X$, on consid\`ere $K(\Omega_{X},2i)$ le complexe de De Rham de 
$X$ ( \cite{h} ), plac\'e en degr\'es $[-2i,DimX-2i]$. 
On prend pour $\cal H\mit^{i}$ un mod\`ele fibrant de $K(\Omega_{-},2i)$.\\
\hspace*{5mm} Notons $Im(X)$ la cat\'egorie filtrante des immersions ferm\'ees
$X \hookrightarrow Y$, avec $Y$ lisse. La cohomologie de $Y$ \`a support dans 
le sous-sch\'ema $X$ est d\'efinie par
$$\cal H\mit^{i}_{X}(Y)=Fib \left(\cal H\mit^{i}(Y) \longrightarrow \cal 
H\mit^{i}(Y-X) \right)$$
L'homologie d'un sch\'ema $X$ est alors d\'efinie par
$$\cal H\mit'_{i}=Colim_{Im(X)}E\cal H\mit^{i+2p}_{X}(Y)$$
o\`u $p$ est la codimension de $X$ dans $Y$, 
et $E\cal H\mit^{i}$ le spectre associ\'e \`a la r\'esolution canonique de 
$\Omega[2i]$ d\'efinie dans \cite[¤$2$]{h}.\\
\hspace*{5mm} Alors $\cal H$ et $\cal H\mit'$ v\'erifient les conditions pr\'ec\'edentes 
avec $d=2$, au d\'etail pr\`es que les images r\'eciproques 
\'etales en homologie ne sont fonctorielles que dans un sens flexible ( 
\cite{s} ). Elles sont d\'efinies de la fa\c{c}on suivante.\\
\hspace*{5mm} Soit $f : X' \rightarrow X$ un morphisme \'etale.
Soit $p : Z_{\bullet} \rightarrow X$ un hyper-recouvrement propre 
dont toutes les composantes $Z_{n}$ sont lisses. Un tel 
hyper-recouvremnent existe d'apr\`es \cite{jo}. En prenant l'image de 
$Z_{\bullet}$ par le foncteur covariant $E\cal H\mit$, on obtient un 
$\Delta$-diagramme de spectres $E\cal H\mit(Z_{\bullet})$. On note la 
colimite homotopique par $E\cal H\mit(Z/X)$. Alors le morphisme 
naturel
$$p_{*} : E\cal H\mit(Z/X) \longrightarrow \cal H\mit'(X)$$
est une \'equivalence faible. On pose $Z_{\bullet}' \rightarrow X'$ 
l'hyper-recouvrement induit sur $X'$. Alors on a un morphisme induit
sur les $\Delta$-diagrammes
$$f^{*} : E\cal H\mit(Z'_{\bullet}) \longrightarrow E\cal 
H\mit(Z_{\bullet})$$
On conclut en prenant la colimite homotopique.
Cette construction et la formule de projection permettent aussi de d\'efinir 
la structure de $\cal H\mit(X)$-module sur $\cal H\mit'(X)$.

\item \underline{Les complexes de Chow-Bloch:} Si l'on se retreint 
aux sch\'emas quasi-projectifs lisses, on peut prendre 
$$\cal H\mit^{i}(X)=\bf Z\mit^{i}(X)\otimes \bf Q$$
le complexe de Bloch des cycles de codimension $i$ sur $X$ ( 
\cite{b} ). En prenant $\cal H\mit'=\cal H$, les conditions 
pr\'ec\'edentes seront v\'erifi\'ees pour les sch\'emas lisses. Si on consid\`ere 
des sch\'emas singuliers, la d\'efinition marche encore sauf que seul
$\cal H\mit'$ existe. On laisse le soin au lecteur d'apporter les 
quelques modifications que cela entraine par la suite. 

\end{itemize}

\begin{df}\label{d3a}
Soit $F$ un champ alg\'ebrique lisse. On d\'efinit la cohomologie et 
l'homologie \'etale 
de $F$ par
$$H_{et}^{p}(F,i):=H^{p-di}(F_{et},\cal H\mit^{i})$$
$$H^{et}_{p}(F,i):=H^{p-di}(F_{et},\cal H\mit'_{i})$$
$$H_{et}^{\bullet}(F,*):=\bigoplus_{i,p}H^{p}_{et}(F,i)$$
$$H^{et}_{\bullet}(F,*):=\bigoplus_{i,p}H_{p}^{et}(F,i)$$
Nous noterons aussi 
$$\cal H\mit(F):=H(F_{et},\cal H\mit)$$
$$\cal H\mit'(F):=H(F_{et},\cal H\mit')$$
les spectres de cohomologie et d'homologie de $F$.
\end{df} 

Nous pr\'ecisons que l'indice $p$ de $H^{et}_{p}(F,*)$ d\'esigne 
la codimension, et non la dimension comme il en est l'usage.\\

Pour d\'efinir les images directes nous aurons besoin du r\'esultat de 
descente analogue \`a \ref{t2b}.

\begin{prop}\label{p3a}
Soit $p : X \rightarrow F$ un morphisme propre et surjectif, avec $X$ 
un sch\'ema quasi-projectif. Alors le morphisme naturel
$$\cal H\mit'(X/F)_{\bf Q} \longrightarrow \cal H\mit'(F)_{\bf Q}$$
est une \'equivalence faible.
\end{prop}

\underline{\bf Preuve:} \rm C'est la m\^eme que pour \ref{t2b}. $\Box$\\

Soit $f : F' \longrightarrow F$ un morphisme propre de champs 
alg\'ebriques. Consid\`erons un diagramme commutatif
$$\xymatrix{
X' \ar[d]_{q} \ar[dr]^{p} &  \\
F' \ar[r]_{f} & F }$$
o\`u $q$ est un morphisme propre et surjectif.
On dispose alors du diagramme suivant
$$\xymatrix{
\cal H\mit'(X'/F')_{\bf Q} \ar[dr]^{p_{*}} \ar[d]_{q_{*}} &  \\
\cal H\mit'(F')_{\bf Q}  & \cal H\mit'(F)_{\bf Q} }$$
Comme $q_{*}$ est une \'equivalence faible, on peut 
d\'efinir $f_{*} : \cal H\mit'(F')_{\bf Q} \longrightarrow \cal 
H\mit'(F)_{\bf Q}$
par 
$$f_{*}:=p_{*} \circ (q_{*})^{-1}$$
On obtient ainsi un morphisme bien d\'efini
$$p_{*} : H^{et}_{\bullet}(F',*)_{\bf Q} \rightarrow 
H^{et}_{\bullet}(F,*)_{\bf Q}$$
Si $p : F \longrightarrow Spec k$ est le morphisme structural d'un 
champ propre, on notera 
$$\int_{F}^{et}:=p_{*} : H^{et}_{\bullet}(F,*)_{\bf Q} \longrightarrow
H^{et}_{\bullet}(Spec k,*)_{\bf Q}$$
La construction de Gillet des classes de Chern ( \cite[$2.2$]{g} ) donne 
des morphismes de pr\'efaisceaux simpliciaux
$$C_{i} : K_{[0]} \rightarrow \cal H\mit^{i}_{[0]}$$
En passant \`a la cohomologie on obtient des classses de Chern \'etale
$$C_{i}^{et} : K_{m,et}(F) \longrightarrow H^{di-m}_{et}(F,i)$$
v\'erifiant les propri\'et\'es cit\'ees dans \cite[$2.1$]{g}.
Sont associ\'es par les formules habituelles, le carat\`ere de Chern et 
les classes de Todd
$$Ch^{et} : K_{*,et}(F) \longrightarrow H^{\bullet}_{et}(F,*)$$
$$Td^{et} : K_{*,et}(F) \longrightarrow H^{\bullet}_{et}(F,*)$$

Nous poss\`edons de cette fa\c{c}on une th\'eorie cohomologie-homologie \'etale pour 
les champs alg\'ebriques, v\'erifiant toutes les propri\'et\'es habituelles ( 
localisation, homotopie ... ).\\

\underline{\bf Remarque:} \rm Voici un exemple tr\`es simple montrant que cette 
th\'eorie \'etale ne peut pas donner de formule de Hirzebruch-Riemann-Roch.\\
\hspace*{5mm} Soit $F=[Spec k/H]$ avec $H$ un groupe fini ab\'elien par exemple. Dans ce
cas le fibr\'e tangent est trivial, et la transformation de Riemann-Roch
$$\tau_{F} : K_{0}(F) \rightarrow H^{d*}_{et}(F,\Gamma(\bullet))_{\bf Q}$$
est le caract\`ere de Chern. Elle est donc multiplicative. Le corollaire
pr\'ec\'edent implique que
$$p_{*} : H^{d*}_{et}(F,\Gamma(\bullet)) \simeq H^{d*}_{et}(Spec k,\Gamma(\bullet))
_{\bf Q}$$
Supposons que la th\'eorie cohomologique v\'erifie
$$H^{d*}_{et}(Spec k,\Gamma(\bullet))_{\bf C}=\bf C$$ 
C'est le cas pour les exemples pr\'ec\'edents.
Si la formule d'Hirzebruch-Riemann-Roch \'etait
v\'erifi\'ee on aurait un diagramme commutatif
$$\xymatrix{
K_{0}(F)_{\bf C} \ar[r]^-{Ch} \ar[d]_{p_{*}} & 
H^{d*}_{et}(F,\Gamma(\bullet))_{\bf C} \ar[d]^{\wr p_{*}}\\ 
\bf C\mit \ar[r]_{Id} & \bf C }$$
En identifiant $K_{0}(F)$ avec le groupe de Grothendieck des 
repr\'esentations de $H$ dans $k$, on aurait pour tout $k[H]$-module
$V$ de dimension finie
$$Ch(V)=Dim(V^{H})$$
Prenons $V$ de dimension $1$ non triviale. Alors $V^{H}=(0)$. Soit 
$m$ l'ordre de $H$, alors $V^{\otimes m}=1$. On aurait donc
$$Ch(V^{\otimes m})=Ch(V)^{m}=Ch(1)=1 \Rightarrow Ch(V)\not = 0$$
Ce qui est absurde.\\

Pour rem\'edier \`a ce probl\`eme nous proposons une nouvelle d\'efinition de 
la cohomologie d'un champ, que nous appellerons ( voir la remarque 
suivante ) "cohomologie \`a coefficients dans les repr\'esentations". Sa 
d\'efinition est directement inspir\'ee du th\'eor\`eme \ref{t2d}.\\

\begin{df}\label{d3b}
Soit $F$ un champ alg\'ebrique. On d\'efinit la cohomologie et 
l'homologie de $F$ \`a coefficients dans les repr\'esentations par
$$H_{rep}^{p}(F,i):=H_{et}^{p}(I_{F},i)$$
$$H^{rep}_{p}(F,i):=H^{et}_{p}(I_{F},i)$$
On notera de m\^eme 
$$\cal H\mit_{rep}(F):=\cal H\mit(I_{F})$$
$$\cal H\mit'_{rep}(F):=\cal H\mit'(I_{F})$$
$$H_{rep}^{\bullet}(F,*):=\bigoplus_{i,p}H^{p}_{rep}(F,i)$$
$$H^{rep}_{\bullet}(F,*):=\bigoplus_{i,p}H_{p}^{rep}(F,i)$$
Si $f : F \longrightarrow F'$ est un morphisme de champs,
et $If : I_{F} \longrightarrow I_{F'}$ le morphisme induit, on 
d\'efinit 
$$f^{*} : H^{\bullet}_{rep}(F,*)=H^{\bullet}_{et}(I_{F},*) 
\stackrel{If^{*}}{\longrightarrow} 
H^{\bullet}_{et}(I_{F'},*)=H^{\bullet}_{rep}(F',*)$$
De la m\^eme fa\c{c}on si $f : F \longrightarrow F'$ est propre, on pose
$$f_{*} : H_{\bullet}^{rep}(F,*)_{\bf Q}=H_{\bullet}^{et}(I_{F},*)_{\bf Q} 
\stackrel{If_{*}}{\longrightarrow} 
H_{\bullet}^{et}(I_{F'},*)_{\bf Q}=H_{\bullet}^{rep}(F',*)_{\bf Q}$$
Si $p : F \longrightarrow Spec k$ est la projection d'un champ
propre, on notera 
$$\int_{F}=p_{*} : H^{rep}_{\bullet}(F,*)_{\bf Q} \longrightarrow 
H^{rep}_{\bullet}(Spec k,*)_{\bf Q}$$
le morphisme induit.
\end{df}

Les propri\'et\'es de cette cohomologie se d\'eduisent imm\'ediatement de 
celles de la cohomologie \'etale. Comme d'habitude, pour tout $\bf 
Z$-module $A$,  nous noterons 
$A_{\Lambda}$ pour $A\otimes_{\bf Z}\Lambda$.\\

\begin{prop}\label{p3b}
La correspondance $F \mapsto H_{rep}^{\bullet}(F,*)_{\Lambda}$ est un foncteur 
contravariant de la $2$-cat\'egorie des champs alg\'ebriques vers les 
$\Lambda$-alg\'ebres commutatives bi-gradu\'ees. \\
\hspace*{5mm} La correspondance $F \mapsto H^{rep}_{\bullet}(F,*)_{\Lambda}$ est un 
foncteur covariant de la $2$-cat\'egorie des champs alg\'ebriques et 
morphismes propres vers 
les $\Lambda$-modules. C'est aussi un foncteur contravariant pour les 
morphismes \'etales repr\'esentables.
De plus, on a les propri\'et\'es suivantes\\
\begin{enumerate}
\item pour tout champ $F$, $H^{rep}_{\bullet}(F,*)_{\Lambda}$ est un 
$H_{rep}^{\bullet}(F,*)_{\Lambda}$-module bi-gradu\'e. Si $f : F' \longrightarrow 
F$ est un morphisme propre, $x \in H^{rep}_{\bullet}(F,*)_{\Lambda}$ 
et $y \in H_{rep}^{\bullet}(F',*)_{\Lambda}$, alors
$$f_{*}(f^{*}(x).y)=x.f_{*}(y)$$
\item si $F$ est lisse, il existe un isomorphisme
de $H_{rep}^{\bullet}(F,*)_{\Lambda}$-module
$$p_{F} : H_{rep}^{\bullet}(F,*)_{\Lambda} \simeq 
H^{rep}_{\bullet}(F,*)_{\Lambda}$$
compatible avec les images r\'eciproques et les produits.
\item si le carr\'e suivant est cart\'esien
$$\xymatrix{
G' \ar[r]^{q} \ar[d]_{v} & F' \ar[d]^{u} \\
G \ar[r]_{p} & F }$$
avec $p$ propre, et $u$ \'etale repr\'esentable, alors
$$q_{*} \circ v^{*} = u^{*} \circ p_{*}$$
\item si $j : F' \hookrightarrow F$ est une immersion ferm\'ee, et 
$i : U \hookrightarrow F$ l'immersion compl\'ementaire, alors la suite 
$$\xymatrix{
\cal H\mit'_{rep}(F') \ar[r]^{j_{*}} & \cal H\mit'_{rep}(F) \ar[r]^{i^{*}} & \cal 
H\mit'_{rep}(U) }$$ 
est naturellement exacte.
\item si $F$ est un champ lisse, et $p : V \rightarrow F$ est un fibr\'e vectoriel, alors le 
morphisme naturel
$$p^{*} : H^{\bullet}_{rep}(F,*) \longrightarrow 
H^{\bullet}_{rep}(V,*)$$
est un isomorphisme.
\end{enumerate}
\end{prop}

\underline{\bf Remarque:} \rm Supposons que $k=\bf C$, et que l'on 
prenne la cohomologie de De Rham. 
Si $F$ est un champ alg\'ebrique lisse, on peut lui associer un champ 
analytique $F^{an}$. Le th\'eor\`eme de comparaison de Grothendieck ( 
\cite[Thm. $1'$]{gr} ), donne alors un isomorphisme de $\bf C$-alg\`ebres
$$H_{rep}^{\bullet}(F)\simeq H^{\bullet}((I_{F})_{top},\bf C\mit)$$
o\`u $(I_{F})_{top}$ est le site topologique sur $(I_{F})^{an}$, et $\bf 
C$ le faisceau constant de fibre $\bf C$. Soit $\pi : I_{F} 
\rightarrow F$ la projection, et  $p : F \rightarrow M$ la projection sur 
l'espace de 
modules. On pose $R=p_{*}\circ \pi_{*}(\bf C\mit)$, c'est un 
faisceau en $\bf C$-alg\`ebres sur $M$, dont la fibre au point 
g\'eom\'etrique $x \in M$ est isomorphe \`a $\bf C\mit(H_{x})$, la $\bf 
C$-alg\`ebre des fonctions centrales sur le groupe d'isotropie $H_{x}$ 
de $x$. On peut donc \'ecrire
$$H_{rep}^{\bullet}(F)\simeq H^{\bullet}(M_{top},R)$$
Cet isomorphisme explique le nom de "cohomologie \`a coefficients dans 
les repr\'esentations", en gardant \`a l'esprit que les \'el\'ements de $\bf 
C\mit(H_{x})$ s'identifient \`a des repr\'esentations virtuelles de 
$H_{x}$, \`a coefficients complexes.

\begin{df}\label{d3c}
Soit $F$ un champ alg\'ebrique. On d\'efinit le caract\`ere de Chern par la 
composition
$$\xymatrix{
K_{m}(F)_{\Lambda} \ar[r]^-{\phi_{F}} & K_{m,et}(I_{F})_{\Lambda} 
\ar[r]^-{Ch^{et}} & 
H^{2i-m}_{et}(I_{F},i)_{\Lambda}=H^{2i-m}_{rep}(F,i)_{\Lambda} }$$
\end{df}

Avant d'\'enoncer le th\'eor\`eme de Riemann-Roch, il nous reste \`a d\'efinir 
la classe de Todd d'un champ lisse.\\

Soit $F$ un champ alg\'ebrique lisse, et $\pi : I_{F} \rightarrow F$ la 
projection canonique. Notons $\cal N\mit^{\vee}=\Omega_{I_{F}/F}$ le 
fibr\'e conormal de $I_{F}$ relativement \`a $F$, $\lambda^{i} : K_{0}(I_{F}) \rightarrow K_{0}(I_{F})$
la $i$-\`eme $\lambda$-op\'eration ( \cite[$V$, ¤$1$]{fl} ), et 
$$\begin{array}{cccc}
\lambda_{-1} : & K_{0}(I_{F}) & \longrightarrow & K_{0}(I_{F}) \\
               & x            & \mapsto         & \sum_{i}(-1)^{i}\lambda^{i}(x) 
\end{array}$$
En utilisant les notations suivant le th\'eor\`eme \ref{t2d}, notons 
$$\alpha_{F}:= can \circ \rho \left( \lambda_{-1}(\cal N\mit^{\vee}) 
\right) \in K_{0,et}(I_{F})_{\Lambda}$$

\begin{lem}\label{l3a}
L'\'el\'ement $\alpha_{F}$ est inversible dans l'anneau 
$K_{0,et}(I_{F})_{\Lambda}$.
\end{lem}

\underline{\bf Preuve:} \rm Soit $m$ le nombre de composantes connexes 
de $I_{F}$, et 
$$rk : K_{0,et}(I_{F})_{\Lambda} \rightarrow \Lambda^{m}$$
l'application rang. On sait alors qu'un \'el\'ement de 
$K_{0,et}(I_{F})_{\Lambda}$ est inversible si et seulement si
son rang est inversible dans $\Lambda^{m}$. Comme ceci est local sur $M_{et}$, 
on peut supposer que $F$ est un champ quotient. 
Le lemme est alors \cite[$1.8$]{vi}.$\Box$\\

\begin{df}\label{d3d}
Soit $F$ un champ alg\'ebrique lisse. On d\'efinit sa classe de Todd par
$$Td(F):=Ch^{et}(\alpha_{F}^{-1}).Td^{et}(T_{I_{F}})$$
La transformation de Riemann-Roch est d\'efinie par
$$\begin{array}{cccc}
\tau_{F} : & K_{m}(F) & \longrightarrow & H^{\bullet}_{rep}(F,*) \\
           & x        & \mapsto         & Ch(x).Td(F)
\end{array}$$
\end{df}

\end{subsection}

\begin{subsection}{D\'emonstration du Th\'eor\`eme}
\hspace{5mm}

On en arrive au th\'eor\`eme de Grothendieck-Riemann-Roch. Pour pouvoir le d\'emontrer 
nous aurons besoin d'hypoth\`eses suppl\'ementaires sur les champs 
que l'on consid\`ere.

\begin{hyp}\label{h2}
Les champs alg\'ebriques $F$ consid\'er\'es, v\'erifieront les hypoth\`eses 
suivantes\\
\begin{enumerate}
\item l'espace de modules $M$ de $F$ est un sch\'ema quasi-projectif.
\item tout faisceau coh\'erent sur $F$ est quotient d'un faisceau 
localement libre.
\end{enumerate}
\end{hyp}

Remarquons que la seconde hypoth\`ese est v\'erifi\'ee dans tous les 
exemples que l'on rencontre ( champs de modules de courbes, de 
vari\'et\'es ab\'eliennes ... ). En effet ces champs sont construits en 
prenant des quotients de sch\'emas quasi-projectifs par des actions de 
groupes r\'eductifs. Et on sait que ces champs v\'erifient la seconde 
hypoth\`ese ( \cite[$5.6$]{th2} ). De plus, elle entraine que le 
morphisme naturel
$$K_{*}(F) \longrightarrow G_{*}(F)$$
est un isomorphisme quand $F$ est lisse.\\
\hspace*{5mm} Soit $\cal L$ un fibr\'e en droite tr\`es ample sur $M$, et 
$\cal L\mit'=p^{*}\cal L$ le fibr\'e induit sur $F$. Alors $\cal L\mit'$ est ample 
sur $F$. En particulier, les deux hypoth\`eses entrainent que tout 
morphisme propre repr\'esentable est projectif ( \cite[$3.3$]{ko} ).

\begin{prop}\label{p3c}
Soit $F$ un champ alg\'ebrique v\'erifiant \ref{h2}, alors 
il existe une enveloppe de Chow rationnelle 
$$q : F_{0} \longrightarrow F$$
avec $F_{0}$ une gerbe triviale sur un sch\'ema quasi-projectif  
et $q$ un morphisme fini.
\end{prop}

\underline{\bf Preuve:} \rm Soit $X \longrightarrow F$ un morphisme 
fini et surjectif, 
avec $X$ un sch\'ema quasi-projectif normal. Un tel morphisme existe d'apr\`es 
\cite[$2.6$]{vi2}. On consid\`ere le morphisme induit sur les espaces de modules $X 
\rightarrow M$, et $F_{X}=F\times_{M}X \rightarrow X$ le changement de 
base. Alors, par construction, le morphisme $X \rightarrow F$ induit une section de la 
projection $F_{X} \longrightarrow X$
$$s : X \longrightarrow F_{X}$$ 
Notons $F_{X}'$ la normalisation de $F_{X}$. On a alors un morphisme 
fini repr\'esentable
$$F_{X}' \longrightarrow F$$
qui v\'erifie la condition de \ref{d2c} pour les points g\'en\'eriques de $F$. 
Or $F_{X}'$ et $X$ sont normaux, et $F_{X}' \longrightarrow X$ poss\`ede une 
section. On sait alors d'apr\`es \cite[$2.7$]{vi2} que $F_{X}'$ est une gerbe sur 
$X$, qui est triviale car elle poss\`ede une section. \\
\hspace*{5mm} Il existe alors un sous-espace ouvert dense $U \hookrightarrow M$, tel que 
le morphisme induit au-dessus de $U$
$$q_{U} : (F_{X})_{U} \longrightarrow F_{U}$$
poss\`ede une section apr\`es un changement de base fini de $U$. En proc\'edant 
alors par r\'ecurrence sur $Dim F$, la proposition est vraie pour le ferm\'e 
compl\'ementaire r\'eduit $F'=F-F_{U}$. Soit $F_{0}' \longrightarrow F'$ 
un morphisme v\'erifiant les conditions demand\'ees. On pose alors
$$F_{0}=F_{X} \coprod F_{0}'$$
Ainsi, le morphisme $q : F_{0} \longrightarrow F$ v\'erifie bien les 
conditions de la proposition. $\Box$\\

\begin{thm}{( Grotendieck-Riemann-Roch )}\label{GRR}
Soit $\cal CH$ la cat\'egorie des champs $F$, v\'erifiants \ref{h} et \ref{h2}.
Alors, pour tout champ $F$ de $\cal CH$, il existe un morphisme de
$\Lambda$-alg\`ebres
$$\tau_{F} : G_{m}(F)_{\Lambda} \longrightarrow 
H_{\bullet}^{rep}(F,*)_{\Lambda}$$
v\'erifiant
\begin{enumerate}
\item $\tau_{F}$ commute avec les images directes de morphismes 
propres de $\cal CH$.
\item $\tau_{F}$ commute avec les images r\'eciproques de morphismes 
repr\'esentables
\'etales de $\cal CH$.
\item si $F$ est un champ lisse de $\cal CH$, alors 
$$\tau_{F}(x)=Ch(x).Td(F)$$
pour tout $x \in K_{*}(F)$.
\item si $F$ est dans $\cal CH$, $x \in K_{*}(F)$, et $y \in 
G_{*}(F)$, alors
$$\tau_{F}(x.y)=Ch(x).\tau_{F}(y)$$
\item si $F$ est un sch\'ema, alors $\tau_{F}$ est le morphisme d\'efini
dans \cite[$4.1$]{g}.
\end{enumerate}
\end{thm}

\underline{\bf Preuve:} \rm Pour ce qui suit, tous les coefficients seront \'etendus \`a $\Lambda$.
Pour simplifier les notations nous n'\'ecrirons pas l'indice $\Lambda$.\\ 

Commen\c{c}ons par d\'emontrer $(1)$ dans le cas des champs lisses.\\

Soit $f : F \longrightarrow F'$ un morphisme propre de champs lisses de $\cal 
CH$, et $If : I_{F} \longrightarrow I_{F'}$ le morphisme induit. Pour un 
champ $G$, on notera $\pi_{G}$ la projection de $I_{G}$ sur $G$.

\begin{lem}\label{l3b}
Le diagramme suivant commute
$$\xymatrix{
G_{*}(F) \ar[d]_-{f_{*}} \ar[r]^-{\alpha_{F}^{-1}.\phi_{F}} & 
G_{*,et}(I_{F}) \ar[d]^-{If_{*}} \\
G_{*}(F') \ar[r]_-{\alpha_{F'}^{-1}.\phi_{F'}} & G_{*,et}(I_{F}) }$$
\end{lem}

\underline{\bf Preuve:} \rm On proc\`ede en deux \'etapes.\\

\underline{Etape $(a)$:} Le morphisme $f$ est repr\'esentable.\\

D'apr\'es l'hypoth\`ese \ref{h2}, $f$ est en fait un morphisme 
projectif. On choisit alors une factorisation
$$\xymatrix{
F \ar[r]^{j} & P \ar[r]^{p} & F' }$$
o\`u $j$ est une immersion ferm\'ee, et $P$ un fibr\'e projectif associ\'e \`a un 
fibr\'e vectoriel $V$ sur $F'$. Par fonctorialit\'e, il nous suffit donc de 
d\'emontrer le lemme pour $j$ et $p$.\\

On sait que
$$\phi_{P}\circ j_{*}=can \circ \rho \circ \pi_{P}^{*} \circ j_{*}$$
Or, comme $j$ est une immersion ferm\'ee, le diagramme suivant est cart\'esien
$$\xymatrix{
I_{F} \ar[r]^{Ij} \ar[d]_{\pi_{F}} & I_{P} \ar[d]^{\pi_{P}} \\
F \ar[r]_{j} & P }$$
La formule d'exc\`es d'intersection ( \cite[$3.8$]{ko} ), donne 
$$\pi_{P}^{*} \circ j_{*}=Ij_{*} \circ \lambda_{-1}(\cal N\mit).\pi_{F}^{*}$$
o\`u $\cal N\mit$ est le fibr\'e virtuel conormal d'exc\`es du diagramme cart\'esien
pr\'ec\'edent  
$$\cal N\mit=\cal N\mit_{F}^{\vee} - Ij^{*}(\cal N\mit_{P}^{\vee})$$
avec $\cal N\mit_{F}^{\vee}$ ( resp. $\cal N\mit_{P}^{\vee}$ ) le fibr\'e 
conormal de $I_{F}$ relativement \`a $F$ ( resp. de $I_{P}$ relativement \`a 
$P$ ).
Ainsi, 
$$\begin{array}{cc}
\phi_{P}\circ j_{*}& =can \circ \rho \circ
\lambda_{-1}(\cal N\mit_{P}^{\vee}).Ij_{*} \circ 
\lambda_{-1}(\cal N\mit_{F}^{\vee})^{-1}.\pi_{F}^{*}\\
                   & =\alpha_{P}.Ij_{*} \circ \alpha_{F}^{-1}.\phi_{F}
\end{array}$$
Ce qui d\'emontre le lemme pour le morphisme $j$.\\

Pour le cas de la projection $p : P \longrightarrow F'$, d'un fibr\'e 
projectif de rang $r$, la formule du fibr\'e projectif implique qu'il suffit de v\'erifier 
que pour tout entier $m$, avec $0 \leq m \leq r$, on a
$$Ip_{*}(\alpha^{-1}_{P}.\phi_{P}(\cal O\mit_{P}(m))=
\alpha^{-1}_{F}.\phi_{F}(p_{*}\cal O\mit_{P}(m))$$
o\`u $\cal O\mit_{P}(1)$ est le fibr\'e canonique sur $P$, et 
$\cal O\mit_{P}(m)=\cal O\mit_{P}(1)^{\otimes m}$. Cette formule se 
v\'erifie par un cacul direct, mais un peu long.\\

Remarquons que la d\'emonstration pr\'ec\'edente reste valable dans le cas o\`u $F$ est une gerbe 
\'eventuellement singuli\`ere. On dispose ainsi du lemme avec $F$ une gerbe 
quelconque, et $f$ un morphisme 
projectif. \\

\underline{Etape $(b)$:} Cas g\'en\'eral.\\

Remarquons d'abors, que lorsque $F$ et $F'$ sont des gerbes triviales ( 
\'eventuellements singuli\`eres ) on a clairement un diagramme commutatif
$$\xymatrix{
G_{*}(F) \ar[r]^{\phi_{F}} \ar[d]_{f_{*}} & G_{*,et}(I_{F}) \ar[d]_{If_{*}} \\
G_{*}(F') \ar[r]_{\phi_{F'}} & G_{*,et}(I_{F'}) }$$

D'apr\'es \ref{p3c} il existe des enveloppes de Chow rationnelles 
finies
$$q : F_{0} \longrightarrow F$$
$$q' : F'_{0} \longrightarrow F'$$
telle que les champs $F_{0}$ et $F'_{0}$ soient des gerbes triviales sur 
des sch\'emas, et qu'il existe un diagramme homotopiquement commutatif
$$\xymatrix{
F_{0} \ar[d]_{q} \ar[r]^{f_{0}} & F_{0}' \ar[d]^{q'} \\
F \ar[r]_{f} & F' }$$
Soient 
$Iq : I_{F_{0}} \longrightarrow I_{F}$, et $Iq' : I_{F'_{0}} \longrightarrow I_{F'}$ 
les morphismes induits. \\

Consid\`erons le diagramme suivant
$$\xymatrix{
 & G_{*}(F_{0}) \ar[dd]_(.3){\phi_{F_{0}}} \ar[rr]^{(f_{0})_{*}} 
 \ar[dl]^{q_{*}} & & G_{*,et}(F_{0}) \ar[dd]^{\phi_{F_{0}'}} 
 \ar[dl]^{q'_{*}} \\
G_{*}(F) \ar[dd]_{\phi_{F}} \ar[rr]^(.7){f_{*}} & & G_{*}(F') 
\ar[dd]^(.3){\phi_{F'}} \\
& G_{*,et}(I_{F_{0}}) \ar[rr]_(.3){(If_{0})_{*}} \ar[dl]^{Iq_{*}}
 & & G_{*,et}(I_{F'_{0}}) \ar[dl]^{Iq'_{0}}\\
G_{*,et}(I_{F}) \ar[rr]_{If_{*}} & & G_{*,et}(I_{F'}) 
}$$

Si nous montrons que le morphisme $q_{*} : G_{*}(F_{0}) \longrightarrow G_{*}(F)$ est
surjectif, il suffit alors de montrer que toutes les faces, sauf la face frontale, commutent. 
Or, il est clair que les deux faces horizontales commutent. La face du fond commute 
car $F_{0}$ et $F_{0}'$ sont des gerbes triviales. Enfin les faces 
lat\'erales commutent d'apr\`es l'\'etape $(a)$. \\

Par d\'efinition d'une enveloppe de Chow rationnelle, le morphisme 
$$I_{q} : I_{F_{0}} \longrightarrow I_{F}$$
est fini et surjectif. Appliquons l'\'etape $(a)$ au morphisme $q$. On en 
d\'eduit un diagramme commutatif
$$\xymatrix{
G_{*}(F_{0}) \ar[r]^{\phi_{F_{0}}} \ar[d]_{q_{*}} & G_{*,et}(I_{F_{0}}) 
\ar[d]^{Iq_{*}} \\
G_{*}(F) \ar[r]_{\phi_{F}} & G_{*,et}(I_{F}) }$$
Le th\'eor\`eme \ref{t2d} implique que $q_{*}$ est surjectif si et seulement 
si $Iq_{*}$ l'est. Comme $I_{F}$ est lisse, la formule de projection 
implique que $Iq_{*}$ est surjectif si le rang de $Iq_{*}(1)$ est inversible dans 
$G_{et,*}(I_{F})$, ce qui est vrai car $Iq$ est fini et surjectif. $\Box$\\
                    
\begin{lem}\label{l3c}
Soit $f : F \longrightarrow F'$ un morphisme propre de champs lisses de
$\cal CH$. Alors, le diagramme suivant commute
$$\xymatrix{
K_{*,et}(F) \ar[d]_-{f_{*}} \ar[rrr]^-{Ch^{et}(-).Td^{et}(T_{F})} & & & 
H^{\bullet}_{et}(F,*)  \ar[d]^-{f_{*}}\\
K_{*,et}(F') \ar[rrr]^-{Ch^{et}(-).Td^{et}(T_{F'})} & & & 
H^{\bullet}_{et}(F',*) }$$
\end{lem}

\underline{\bf Preuve:} \rm La preuve suit le m\^eme sch\'ema que celle du 
lemme pr\'ec\'edent.\\

\underline{Etape $(a)$:} Le morphisme $f$ est repr\'esentable. \\

On sait qu'il est alors projectif. Dans ce cas la d\'emonstration est presque mot 
pour mot celle de \cite[$4.1$]{g}. En effet les seules propri\'et\'es dont on a 
besoin sont, la d\'eformation vers le cone normal, et la structure de la 
$G$-th\'eorie d'un fibr\'e projectif.\\

\underline{Etape $(b)$:} Cas g\'en\'eral. \\

Soit $p : X \rightarrow F'$ un morphisme propre et g\'en\'eriquement 
fini, avec $X$ un sch\'ema lisse. Un tel morphisme $p$ existe d'apr\'es 
\cite[$4.12$]{dm} et \cite{jo}. Alors les morphismes
$$p_{*} : G_{*,et}(X) \longrightarrow G_{*,et}(F)$$
sont surjectifs. De plus $f \circ p : X \rightarrow F'$ est repr\'esentable. 
En utilisant l'\'etape $(a)$ pour $p$ et $f \circ p$, on termine la preuve du 
lemme. $\Box$\\

Revenons \`a notre morphisme $f : F \longrightarrow F'$. Par d\'efinition, 
$$\tau_{F}(-)=Ch^{et}(\alpha_{F}^{-1}.\phi_{F}(-)).Td^{et}(T_{I_{F}})$$
Ainsi le lemme \ref{l3c} appliqu\'e \`a $If : I_{F} \longrightarrow I_{F'}$, donne
$$f_{*}(\tau_{F}(-))=Ch^{et}(f_{*}(\alpha_{F}^{-1}.\phi_{F}(-))).Td^{et}(T_{I_{F'}})$$
Une application du lemme \ref{l3b} au second membre donne
$$\begin{array}{cl}
f_{*}(\tau_{F}(-)) 
& =Ch^{et}(\alpha_{F'}^{-1}.\phi_{F'}(f_{*}(-))).Td^{et}(T_{I_{F'}}) \\ 
& =Ch^{et}(\alpha_{F'})^{-1}.Td^{et}(T_{I_{F'}}).Ch^{et}(\phi_{F}(f_{*}(-))) \\
& =\tau_{F'}(f_{*}(-)) 
\end{array}$$

Indiquons bri\`evement comment on \'etend $\tau$ aux champs 
singuliers. \\

Si $F$ est un champ de $\cal CH$, on va construire $\tau_{F}$ comme le 
compos\'e
$$G(F) \stackrel{\psi_{F}}{\longrightarrow} G_{et}(I_{F})
\stackrel{\tau^{et}_{I_{F}}}{\longrightarrow} H^{rep}_{\bullet}(F,*)$$
Soit $p : F_{\bullet} \longrightarrow F$ une hyper-enveloppe de Chow 
rationnelle, telle que chaque $F_{m}$ soit une gerbe triviale 
$X_{m}\times BH_{m}$. Chaque morphisme $F_{n} \rightarrow F_{m}$ 
d\'efini deux morphismes
$$X_{n} \rightarrow X_{m}$$
$$H_{n} \rightarrow H_{m}$$
et donc des morphismes d'inductions sur les alg\`ebres de fonctions 
centrales \`a valeurs dans $\Lambda$
$$\Lambda(H_{n}) \rightarrow \Lambda(H_{m})$$
On peut donc d\'efinir un spectre simplicial
$$G(R) : [n] \rightarrow G(X_{n})\otimes \Lambda(H_{n})$$
Notons que l'on a une \'equivalence faible canonique
$$\psi_{F_{n}} : G(F_{n})_{\Lambda} \longrightarrow G(X_{n})\times 
\Lambda(H_{n})$$ 
qui est compatible avec les images directes, et donc une \'equivalence de 
spectres simpliciaux
$$\psi_{F_{\bullet}} : G(F_{\bullet}/F)_{\Lambda} \longrightarrow G(R)$$
D'autre part, les espaces de modules $MI_{F_{n}}$ de $I_{F_{n}}$ sont naturellement 
isomorphes \`a $X_{n}\times c(H_{n})$, o\`u $c(H_{n})$ est l'ensemble des 
classes de conjugaisons du groupe $H_{n}$. Ainsi on a une \'equivalence 
faible naturelle
$$G(MI_{F_{n}})_{\Lambda} \simeq G(X_{n})\otimes \Lambda(H_{n})$$
De plus, les morphismes $I_{F_{n}} \longrightarrow F$ d\'efinissent des 
morphismes $MI_{F_{n}} \longrightarrow M$, et donc un morphisme de 
spectre simplicial
$$q_{*} : G(R) \longrightarrow G(M)_{\Lambda}$$ 

Le diagramme suivant
$$\xymatrix{
G(F_{\bullet}/F)_{\Lambda} \ar[d]_{p_{*}} \ar[r]^-{\phi_{F_{\bullet}}} & 
G(R) \ar[d]^{q_{*}} \\
G(F)_{\Lambda} & G(MI_{F})_{\Lambda} }$$
et la proposition \ref{p2d} appliqu\'ee \`a $p$, 
permettent alors de poser
$$\psi_{F}=q_{*} \circ \phi_{F_{\bullet}} \circ (p_{*})^{-1} : G_{*}(F)_{\Lambda} \longrightarrow 
G_{*}(MI_{F})_{\Lambda}$$ 
Si $r : I_{F} \longrightarrow I_{M}$ est la projection 
de $I_{F}$ sur son espace de modules, on d\'efinit le second par le 
diagramme commutatif
$$\xymatrix{
G_{*,et}(I_{F})_{\Lambda} \ar[r]^{\tau^{et}_{I_{F}}} \ar[d]_{\wr r_{*}} & 
H^{rep}_{\bullet}(F,*)_{\Lambda} \ar[d]^{\wr r_{*}} \\ 
G_{*}(IM)_{\Lambda} \ar[r]_{\tau_{IM}} & H_{\bullet}(IM,*)_{\Lambda} }$$
o\`u $\tau_{IM}$ est d\'efini dans \cite[$4.1$]{g}.\\

On v\'erifie ais\'ement que $\psi_{F}$ et $\tau_{I_{F}}^{et}$ sont bien 
d\'efinis, et que $\tau_{F}$ poss\`ede les propri\'et\'es $(2)$, $(3)$ et $(4)$ 
du th\'eor\`eme. 
$\Box$\\

\end{subsection}

\begin{subsection}{Exemples}
\hspace{5mm}
Donnons quelques exemples d'applications du th\'eor\`eme 
\ref{GRR}.\\

\begin{cor}\label{HRR}{( Hirzebruch-Riemann-Roch )}
Soit $\cal F$ un faisceau coh\'erent sur un champ propre et lisse $F$ de 
$\cal CH$. Alors on a
$$\chi(F,\cal F\mit)=\int_{F}Ch(\cal F\mit).Td(F)$$
En particulier
$$\chi(M,\cal O\mit_{M})=\int_{F}Td(F)$$
\end{cor}

\underline{\bf Preuve:} \rm On fait $F'=Spec k$ dans \ref{GRR} $(1)$. 
$\Box$\\ 

Comme application de cette formule, on donne plusieurs formules de type 
Gauss-Bonnet, calculant des carat\'eristiques d'Euler.\\

Supposons que $F$ soit un champ de $\cal CH$, lisse et propre sur $Spec \bf C$. On applique 
la formule d'Hirzebruch-Riemann-Roch au complexe de De Rham de $F$. 
Pour cela, on pose pour un fibr\'e en droite $\cal L$
$$C_{1}(\cal L\mit).Td(\cal L\mit)^{-1}:=1-Ch(\cal L\mit^{-1})$$
Par le "splitting principle", on peut \'etendre cette d\'efinition \`a
tout fibr\'e vectoriel $V$, de fa\c{c}on \`a ce que si 
$V=\bigoplus_{i}\cal L\mit_{i}$, on ait
$$C_{max}(V).Td(V)^{-1}=\prod_{i}C_{1}(\cal L\mit_{i}).Td(\cal 
L\mit_{i})^{-1}$$
Remarquons alors que l'on a la formule suivante
$$Ch(\lambda_{-1}(V^{\vee}))=C_{max}(V).Td(V)^{-1}$$

\begin{cor}
Soit $M$ l'espace de modules d'un champ alg\'ebrique $F$ complexe, lisse et propre
de $\cal CH$. Alors
$$\chi(M_{top})=\int_{F}C_{max}(T_{F}).Td(T_{F})^{-1}.Td(F)$$
\end{cor} 

Dans le cas o\`u $F$ est un champ alg\'ebrique complexe de dimension $1$, propre, 
lisse et g\'en\'eriquement non-ramifi\'e, on 
peut expliciter la formule de Riemann-Roch, et retrouver la formule 
de Gauss-Bonnet d\'emontr\'ee dans \cite[$3.2.5.2$]{ta}.

\begin{cor}
Si $F$ est un champ alg\'ebrique de dimension $1$, lisse et propre sur $Spec \bf C$, appartennant \`a 
$\cal CH$, et tel que $F$ soit g\'en\'eriquement un sch\'ema. On note $g$ le 
genre de l'espace de modules $M$, 
et $x_{1}, \dots x_{r}$ les points de ramification de $F$, d'ordres respectifs 
$n_{1}, \dots n_{r}$.
Alors on a
$$\int^{et}_{F}C^{et}_{1}(T_{F})=2-2g-\sum_{j=1}^{i=r}\frac{n_{j}-1}{n_{j}}$$
\end{cor}

\underline{\bf Preuve:} \rm Soit $\zeta_{j}=e^{\frac{2i\pi}{n_{j}}} \in 
\bf C$. Notons $f_{j} : \widetilde{x_{j}} \rightarrow F$ le morphisme naturel. 
Alors $f_{j}^{*}(T_{F})$ est un fibr\'e vectoriel sur 
$\widetilde{x_{j}}\simeq B(\bf Z\mit/n_{j})$. Il s'identifie \`a la 
repr\'esentation de dimension $1$ de $\bf Z\mit/n_{j}$, qui \`a $k$ fait 
correspondre la mutiplication par $\zeta_{j}^{k}$. Ainsi, la formule 
de Hirzebruch-Riemann-Roch appliqu\'ee \`a $\cal O\mit_{F}$, donne
$$\frac{1}{2}\int_{F}^{et}C^{et}_{1}(T_{F})=\chi(F,\cal 
O\mit_{F})-\sum_{j=1}^{j=r}\frac{e_{j}}{n_{j}}$$
o\`u $e_{j}=\sum_{k=1}^{k=n_{j}-1}\frac{1}{1-\zeta^{-k}}$. Or, si on 
note $P_{j}(X)=1+X+ \dots +X^{n_{j}-1}$, on a 
$$e_{j}=\frac{P'(1)}{P(1)}=\frac{n_{j}(n_{j}-1)}{2n_{j}}=\frac{n_{j}-1}{2}$$
De plus, $\chi(F,\cal O\mit_{F})=\chi(M,\cal O\mit_{M})=1-g$. $\Box$\\
\end{subsection}

Plus g\'en\'eralement, le lemme \ref{l3c} appliqu\'e \`a la projection
$$p : F \longrightarrow Spec \bf C$$
d'un champ complexe lisse et propre de $\cal CH$,  et au complexe de De 
Rham de $F$, permet de retrouver la formule de "Gauss-Bonnet" ( 
\cite[$3.1.4.2$]{ta} ) g\'en\'erale.

\begin{cor}
Soit $F$ un champ alg\'ebrique complexe, lisse et propre de $\cal CH$. Notons
$$M_{d} \hookrightarrow M_{d-1} \hookrightarrow \dots \hookrightarrow 
M_{0}=M$$
une stratification de l'espace de modules, telle que $F$ soit une gerbe 
$F_{i}$ sur 
$\widetilde{M_{i}}:=M_{i}-M_{i+1}$, avec $F_{i}$ irr\'eductible. Alors
$$\int_{F}^{et}C_{n}^{et}(T_{F})=
\chi^{orb}(F):=\sum_{i}\frac{\chi((\widetilde{M_{i}})_{top},\bf C\mit
)}{n_{i}}$$
o\`u $n_{i}$ est l'ordre d'inertie de la gerbe $F_{i}$, et $n=Dim_{\bf C}F$.
\end{cor}

Il existe encore une autre caract\'eristique d'Euler, "la caract\'eristique 
d'Euler des physiciens". On peut la d\'efinir pour un champ 
alg\'ebrique complexe $F$ par
$$\chi^{phy}(F):=\sum_{i}(-1)^{i}Dim_{\bf C}H^{i}_{rep}(F)$$
o\`u l'on utilise la cohomologie de De Rham. Cette d\'efinition est 
compatible avec la d\'efinition donn\'ee dans \cite{as} pour un 
quotient par un groupe fini.\\

Pour finir remarquons que le morphisme $\psi$ d\'efini \`a la fin de la 
preuve du th\'eor\`eme, reste valable si l'on suppose simplement que 
$F$ v\'erifie \ref{h}. On dispose alors d'un morphisme naturel
$$\psi_{F} : G_{*}(F) \longrightarrow G_{*}(MI_{F})$$
Soit $r : I_{F} \longrightarrow MI_{F}$ la projection. On d\'efinit $\tau_{I_{F}}^{et}$ par 
le carr\'e commutatif suivant
$$\xymatrix{
G_{0,et}(I_{F})_{\Lambda} \ar[r]^{\tau^{et}_{I_{F}}} \ar[d]_{\wr r_{*}} & 
H^{rep}_{\bullet}(F,*)_{\Lambda} \ar[d]^{\wr r_{*}} \\ 
G_{0}(I_{M})_{\Lambda} \ar[r]_{\tau_{I_{M}}} & H_{\bullet}(I_{M},*)_{\Lambda} }$$
o\`u $\tau_{I_{M}}$ est d\'efini dans \cite[$8.3$]{g3}. 
Ainsi, par composition on a une transformation de Riemann-Roch
$$\tau_{F}=\tau_{I_{F}}^{et} \circ \psi_{F} : G_{0}(F)_{\Lambda} \longrightarrow 
H^{rep}_{\bullet}(F,*)_{\Lambda}$$
qui \'etend celle du th\'eor\`eme \ref{GRR}, au cas des champs alg\'ebriques 
ne v\'erifiant que \ref{h}.\\

\end{section}
\begin{section}{Appendice}\label{app}
\hspace{5mm}

Dans cet appendice on d\'emontre le th\'eor\`eme de descente. 
La preuve est d\'ej\`a dans la preuve de \cite[Thm. $3-10$]{j2}, mais 
le r\'esultat n'\'etant pas explicit\'e sous cette forme nous avons tenu
\`a en donner une d\'emonstration compl\`ete.\\

Pour tout l'appendice, $C$ est un site. On utilisera les notations de
la section $1$, ainsi que la proposition suivante :

\begin{prop}{\cite[$I-1$ Cor. $1$]{q2}}\label{a1}
Soit $F$ un pr\'efaisceau simplicial fibrant, et $f : A \rightarrow B$
une \'equivalence faible. Alors le morphisme induit
$$f^{*} : Hom_{s}(B,F) \rightarrow Hom_{s}(A,F)$$
est une \'equivalence faible.
\end{prop}
 
Si $X\in C$, alors on dispose d'un foncteur image r\'eciproque
$$j^{*} : SPr(C) \rightarrow SPr(C/X)$$
Ce foncteur poss\`ede un adjoint \`a gauche 
$$j_{!} : SPr(C/X) \rightarrow SPr(C)$$
qui est l'extension par le pr\'efaisceau vide. Il est d\'efini 
par :\\
\hspace*{5mm} pour $F\in SPr(C/X)$ et $U\in C$, alors
$$(j_{!}F(U)):=\coprod_{Hom_{C}(U,X)}F(U\rightarrow X)$$
Il est clair que $j_{!}$ pr\'eserve les cofibrations ainsi que 
les \'equivalences faibles.

\begin{lem}\label{a2}
Soit $C$ et $C'$ deux sites et un foncteur
$$a : SPr(C') \rightarrow SPr(C)$$ 
poss\`edant un adjoint \`a gauche
$$b : SPr(C) \rightarrow SPr(C')$$
qui pr\'eserve les cofibrations et les \'equivalences faibles, et tel que
$b(*)=*$. Alors le foncteur $a$ transforme objets fibrants en objets 
fibrants.
\end{lem}

\underline{\bf Preuve:} \rm Soit $F$ un pr\'efaisceau simplicial fibrant sur 
$C'$, et $i : A \hookrightarrow B$ une cofibration triviale de $SPr(C)$. Il 
faut montrer que le morphisme induit
$$i^{*} : Hom(B,a(F)) \rightarrow Hom(A,a(F))$$
est surjectif. Mais par adjonction, on dispose d'un carr\'e commutatif
$$\begin{array}{ccc}
Hom(B,a(F)) & \rightarrow & Hom(A,a(F)) \\
 \wr \downarrow & & \downarrow \wr\\
Hom(b(B),F) & \rightarrow & Hom(b(B),F)
\end{array}$$
Or par hypoth\`ese $b(i) : b(A) \rightarrow b(B)$ est une cofibration 
triviale de $SPr(C')$, et donc le morphisme du bas est surjectif. Ce qui 
implique que celui du haut aussi. $\Box$\\
 
On vient de voir que si $F$ est fibrant sur $C$, alors le 
pr\'efaisceau $j^{*}F$ ( que l'on notera $F_{X}$  par la suite ) est fibrant 
sur $C/X$. De plus, comme $j^{*}$ pr\'eserve les cofibrations triviales, 
on obtient une \'equivalence faible canonique
$$H(C/X,F_{X}) \stackrel{\simeq}{\rightarrow} F^{\circ}(X)$$
pour chaque r\'esolution injective $F \hookrightarrow F^{\circ}$. De cette 
fa\c{c}on on identifiera toujours les espaces $F^{\circ}(X)$ et 
$H(C/X,F_{X})$, que l'on notera $H(X,F)$.\\ 

\begin{df}\label{a3}
Un objet simplicial $X_{\bullet}$ de $C$ est un foncteur 
$$X_{\bullet} : \Delta^{op} \rightarrow C$$
o\`u $\Delta$ est la cat\'egorie simpliciale standart. 
On notera $X_{m}$ pour l'objet $X_{\bullet}([m])$.\\
\hspace*{5mm} Si $X_{\bullet}$ est un objet simplicial de $C$, le site induit sur 
$X_{\bullet}$ est le site suivant :
\begin{itemize}
\item les objets sont les morphismes de $C$
$$U \rightarrow X_{m}$$
pour $m$ un entier positif.
\item un morphisme de $f : U \rightarrow X_{m}$ vers $g : V \rightarrow 
X_{n}$ est la donn\'ee d'un morphisme $a : [n] \rightarrow [m]$ dans 
$\Delta$ et d'un diagramme commutatif dans $C$
$$\begin{array}{rcl}
U & \rightarrow & V \\
f \downarrow & & \downarrow g \\
X_{m} & \stackrel{X_{\bullet}(a)}{\rightarrow} & X_{n}
\end{array}$$
\item un morphisme est couvrant si le morphisme induit 
$$U \rightarrow V$$ 
est couvrant dans $C$.
\end{itemize}
Ce site est not\'e $C/X_{\bullet}$
\end{df}

Remarquons que l'on a un foncteur de restriction
$$j^{*} : SPr(C) \rightarrow SPr(C/X_{\bullet})$$
A travers ce foncteur, tout pr\'efaisceau simplicial $F$ sera aussi 
consid\'er\'e comme pr\'efaisceau sur $X_{\bullet}$. Ainsi 
$H(X_{\bullet},F)$ d\'esignera $H(C/X_{\bullet},j^{*}F)$.\\

Soit $X$ un objet de $C$ et $U \rightarrow X$ un morphisme couvrant. Le nerf 
du recouvrement $U/X$ est l'objet simplicial de $C$ d\'efini par
$$\begin{array}{ccc}
\Delta & \rightarrow & C \\
 \ 
[m] & \mapsto & 
U^{(m)}=\underbrace{U_{\stackrel{\times}{X}}U \dots 
_{\stackrel{\times}{X}}U}_{m+1 \; fois}
\end{array}$$
et les morphismes $U^{(m)} \rightarrow U^{(n)}$ sont induits par les 
projections et les diagonales. On le note $\cal N\mit(U/X)$.\\

Si $F$ est un pr\'efaisceau simplicial et $U \rightarrow X$ un recouvrement 
de $C$, on obtient un $\Delta$-diagramme ( espace cosimplicial ) 
d'ensembles simpliciaux 
$$\begin{array}{ccc}
\Delta & \rightarrow & SEns \\
 \ [m] & \mapsto & F(U^{(m)})
\end{array}$$

Rappelons que l'espace de cohomologie de $\Check{C}ech$ du recouvrement 
$U/X$ \`a coefficients dans le pr\'efaisceau simplicial $F$ est 
$$\Check{H}(U/X,F):=Holim_{[m] \in \Delta}F(U^{(m)})$$

\begin{df}\label{a4}
Un espace cosimplicial $Z$ est un foncteur
$$\begin{array}{cccc}
Z : & \Delta & \rightarrow & SEns \\
              & [m]   & \mapsto & Z([m])
\end{array}$$
La cat\'egorie des espaces cosimpliciaux est not\'ee $CSEns$ 
\end{df}

\underline{\bf Exemples}\rm 
\begin{itemize}
\item $*$ est l'espace cosimplicial constant
\item l'espace cosimplicial $\Delta/-$ est d\'efini par 
$$[m] \mapsto (\Delta/-)([m])=B(\Delta/[m])$$
o\`u $B(I)$ est l'espace classifiant de la cat\'egorie $I$, 
et $I/i$ est la cat\'egorie des morphismes de $I$ de but $i$
\end{itemize}

Un espace cosimplicial $Z$ peut \^etre vu comme un 
pr\'efaisceau simplicial sur $\Delta$ ( site trivial ). Ainsi, si 
$Y$ et $Z$ sont deux espaces cosimpliciaux, on d\'efinit l'espace des 
morphismes de $Y$ vers $Z$ par 
$$Hom_{cs}(Y,Z):=Hom_{s}(Y,Z)$$
o\`u $Y$ et $Z$ sont consid\'er\'es comme pr\'efaisceaux sur $\Delta$.\\

Avec ces notations, la limite homotopique de $Z$ est donn\'ee par 
( \cite[Ch. $XI$ $3-2$]{bk} )
$$Holim_{\Delta}Z=Hom_{cs}(\Delta/-,Z)$$

\begin{thm}\label{a5}
Si $F$ est un pr\'efaisceau simplicial sur $C$, et $X_{\bullet}$ un objet 
simplicial de $C$. Alors il existe une \'equivalence faible fonctorielle en 
$F$
$$H(X_{\bullet},F) \simeq Holim_{[m]\in \Delta}H(X_{m},F)$$ 
\end{thm}

\begin{lem}\label{a6}
Si $F$ est un objet fibrant de $SPr(C)$, alors le pr\'efaisceau induit $F$ 
sur le site $C/X_{\bullet}$ est flasque.
\end{lem}

\underline{\bf Preuve:} \rm Soit $j^{*} : SPr(C) \rightarrow 
SPr(C/X_{\bullet})$ le morphisme de restriction. On consid\`ere $j^{*}F 
\rightarrow F^{\circ}$ une r\'esolution injective sur $C/X_{\bullet}$.
Alors, d'apr\`es le lemme \ref{a2}, le morphisme induit sur $C/X_{m}$
$$F \rightarrow F^{\circ}$$
est une cofibration triviale d'objets fibrants. C'est donc une \'equivalence 
d'homotopie. Ainsi, pour tout objet $U$ de $C/X_{m}$, le morphisme induit
$$F(U) \rightarrow F^{\circ}(U)$$
est une \'equivalence faible. Comme ceci est vrai pour tout $U$ et tout $m$, 
on en d\'eduit que pour tout objet $U$ de $C/X_{\bullet}$
$$F(U) \rightarrow F^{\circ}(U)$$
est une \'equivalence faible. $\Box$\\

\begin{lem}\label{a7}
Le foncteur 
$$\begin{array}{ccc}
SPr(C/X_{\bullet}) & \rightarrow & CSEns \\
F & \mapsto & F(X_{\bullet})
\end{array}$$
admet un adjoint \`a gauche not\'e $Z \mapsto \widetilde{Z}$
\end{lem}

\underline{\bf Preuve:} \rm Soit $Z$ un espace cosimplicial. On d\'efinit 
$$\begin{array}{cccc}
\widetilde{Z} : & C/X_{\bullet} & \rightarrow & SEns \\
                & ( U \rightarrow X_{m} ) & \mapsto & Z([m])
\end{array}$$
$\Box$ \\

\underline{\bf Preuve du th\'eor\`eme:} \rm Soit $F$ un pr\'efaisceau 
simplicial sur $C$. En rempla\c{c}ant $F$ par $F^{\circ}$ on peut supposer 
que $F$ est fibrant sur $C$.
Soit $F \hookrightarrow F^{\circ}$ une r\'esolution injective sur 
$C/X_{\bullet}$. On sait que $H(X_{\bullet},F)=Hom_{s}(*,F^{\circ})$. 
Notons $\bf1\mit=\widetilde{(\Delta/-)}$. \\

\begin{lem}\label{a8}
Le morphisme canonique $\bf 1\mit \rightarrow *$ est une \'equivalence faible
dans $SPr(C/X_{\bullet})$.
\end{lem}

\underline{\bf Preuve:} \rm Il suffit de voir que pour chaque $m$ l'espace 
$\bf 1\mit(X_{m})$ est contractile. Or par d\'efinition, on a 
$$\bf 1\mit(X_{m})=B(\Delta/[m])$$
Mais le classifiant d'une cat\'egorie qui poss\`ede un objet final est 
contractile. $\Box$ \\

Comme $F^{\circ}$ est fibrant, la proposition \ref{a1} montre que le morphisme 
naturel
$$Hom_{s}(*,F^{\circ}) \stackrel{\simeq}{\rightarrow} Hom_{s}(\bf 
1\mit,F^{\circ})$$
est une \'equivalence faible.
Mais, par adjonction, il existe une \'equivalence faible 
fonctorielle
$$Hom_{s}(\bf 1\mit,F^{\circ}) \stackrel{\simeq}{\rightarrow} 
Holim_{\Delta}F^{\circ}(X_{\bullet})$$
On obtient ainsi une \'equivalence naturelle
$$H(X_{\bullet},F) 
\stackrel{\simeq}{\rightarrow}Holim_{\Delta}F^{\circ}(X_{\bullet})$$
De plus, par le lemme \ref{a6}, le morphisme naturel 
$$F \rightarrow F^{\circ}$$  
est une \'equivalence faible \em objet par objet\rm . Comme
les limites homotopiques pr\'eservent les \'equivalences faibles, le 
morphisme induit
$$Holim_{\Delta}F(X_{\bullet}) \rightarrow 
Holim_{\Delta}F^{\circ}(X_{\bullet})$$
est une \'equivalence faible. Enfin, comme $F$ est fibrant, et par la 
remarque suivant \ref{a2}, on obtient un diagramme fonctoriel en $F$
$$H(X_{\bullet},F) \stackrel{\simeq}{\rightarrow} 
Holim_{\Delta}F^{\circ}(X_{\bullet}) \stackrel{\simeq}{\leftarrow}
Holim_{\Delta}H(X_{m},F)$$
$\Box$\\

\begin{prop}\label{a9}
Si $F$ est un pr\'efaisceau simplicial, alors le morphisme naturel
$$H(X,F) \rightarrow H(\cal N\mit(U/X),F)$$
est une \'equivalence faible
\end{prop}

Pour cela on a besoin d'un lemme que nous ne d\'emontrerons pas ( voir
\cite[Cor. $2-7$]{j2} ).

\begin{lem}\label{a10}
Si $F$ est un faisceau simplicial sur un site $C$, alors il existe une 
r\'esolution injective $F \hookrightarrow F^{\circ}$, o\`u $F^{\circ}$ est 
un faisceau d'ensembles simpliciaux.
\end{lem}

\underline{\bf Preuve de la proposition:} \rm Remarquons que, si l'on note 
$SS(X)$ la cat\'egorie des faisceaux simpliciaux sur le site $C/X$, alors
on dipose d'une \'equivalence de cat\'egories
$$b : SS(X) \rightarrow SS(\cal N\mit(U/X))$$
De plus, si $F \in SS(X)$, qui est aussi fibrant 
comme pr\'efaisceau simplicial, $b(F)$ est alors fibrant en tant qu'objet 
de $SPr(\cal N\mit(U/X))$. En effet, notons $a$ le foncteur 
de faisceautisation.\\
\hspace*{5mm} Soit $A \hookrightarrow B$ une cofibration triviale de $SPr(\cal 
N\mit(U/X))$, et un diagramme commutatif sur $\cal N\mit(U/X)$
$$\begin{array}{ccc}
A & \rightarrow & b(F) \\
\downarrow & & \downarrow \\
B & \rightarrow & *
\end{array}$$
Comme $F$ est un faisceau on obtient un carr\'e commutatif
$$\begin{array}{ccc}
Hom_{\cal N\mit(U/X)}(B,b(F)) & \rightarrow & Hom_{\cal N\mit(U/X)}(A,b(F)) 
\\\| & & \| \\
Hom_{\cal N\mit(U/X)}(a(B),b(F)) & \rightarrow & Hom_{\cal 
N\mit(U/X)}(a(A),b(F))
\end{array}$$
Puis, par l'\'equivalence de cat\'egories $SS(X) \rightarrow SS(\cal 
N\mit(U/X))$, un diagramme commutatif
$$\begin{array}{ccc}
Hom_{\cal N\mit(U/X)}(B,b(F)) & \rightarrow & Hom_{\cal N\mit(U/X)}(A,b(F)) 
\\
\| & & \| \\
Hom_{X}(a(B),F) & \rightarrow & Hom_{X}(a(A),F) \\
\| & & \| \\
Hom_{X}(B,F) & \rightarrow & Hom_{X}(A,F)
\end{array}$$
et le morphisme du bas est surjectif par hypoth\`ese.\\

Soit $F$ un pr\'efaisceau simplicial sur $C$. Quitte \`a remplacer $F$ par 
son faisceau associ\'e, on peut supposer que $F$ est un faisceau. En effet 
le morphisme naturel $F \rightarrow a(F)$ est une \'equivalence faible, donc
$F$ et $a(F)$ ont la m\^eme cohomologie.\\

Soit $F \hookrightarrow F^{\circ}$ une r\'esolution injective dans
$SPr(C)$, avec 
$F^{\circ}$ un faisceau. Alors 
$$b(F) \rightarrow b(F^{\circ})$$
est encore une r\'esolution injective. On en d\'eduit donc 
$$H(\cal N\mit(U/X),F)=Hom_{s}(*,b(F^{\circ}))=Hom_{s}(*,F^{\circ})
=H(X,F)$$
$\Box$\\

\begin{cor}\label{desc}
Si $U \rightarrow X$ est un recouvrement de $C$ et $F$ un pr\'efaisceau 
simplicial, alors il existe une \'equivalence faible fonctorielle
en $F$, $X$ et $U$
$$H(X,F) \stackrel{\simeq}{\rightarrow} \Check{H}(U/X,H(-,F))$$
\end{cor} 

\underline{\bf Preuve:} \rm On applique la proposition \ref{a9} et le 
th\'eor\`eme \ref{a5}. $\Box$\\

\end{section}

\end{document}